\documentclass{amsart}

\usepackage{amsmath,amsthm,amssymb}
\usepackage{enumerate}
\usepackage{graphicx}
\usepackage{xspace}
\usepackage{mathrsfs}
\usepackage{diagrams}
\usepackage[numbers]{natbib}
\usepackage[pdfauthor={Erik Stokes},
            pdftitle={The $h$-vectors of 1-dimensional Matroids Complexes},
            pdfsubject={commutative algebra},
            pdfkeywords={simplicial complex, matroid, $h$-vector,
              arithmetic degree, Stanley-Reisner ideal},
            pdfproducer={Latex with hyperref},
            pdfcreator={latex->dvips->ps2pdf},
            pdfpagemode=UseOutlines,
            bookmarksopen=true,
            bookmarksnumbered=true]{hyperref}
\usepackage{rotating}%provides sidewaystable for sides-ways tables
\newcommand{\subs}{\subseteq} %rename \subseteq to \subs
\newcommand{\abs}[2][ ]{\lvert #2\rvert_{#1}}

\newcommand{\isom}{\cong}%Isomorphism
\renewcommand{\vec}[1]{\mathbf{#1}} %Changes vectors from arrows to bold
\newcommand{\deq}{\mathrel{\mathop:}=}%defined equal to
\newcommand{\N}{\ensuremath{\mathbb{N}}}%the naturals

\newcommand{\ideal}[1]{\left\langle#1\right\rangle}%ideal generated by #1
\newcommand{\defword}[1]{\textit{#1}}%the word being defined
\newcommand{\m}{\mathfrak{m}}
\newcommand{\note}[2][ ]{ }%dummy macro

\newcommand{\restrict}[2]{{#1}\vert_{#2}}
\newcommand{\skeleton}[2]{\left[#2\right]_{#1}}
\DeclareMathOperator{\link}{link}

\newtheorem{thrm}{Theorem}[section]
\newtheorem{lemma}[thrm]{Lemma}
\newtheorem{prop}[thrm]{Proposition}
\newtheorem{cor}[thrm]{Corollary}
\newtheorem*{conj*}{Conjecture}
\newtheorem{conj}[thrm]{Conjecture}
\theoremstyle{definition}
\newtheorem{defn}[thrm]{Definition}
\newtheorem*{notation}{Notation}
\newtheorem{rmrk}[thrm]{Remark}
\newtheorem{ex}[thrm]{Example}
\begin{document}
%
%\commentmode
%

%
\title[1-dimensional Matroid Complexes]{The $h$-vectors of
  1-dimensional Matroid Complexes and a 
  Conjecture of Stanley}
\author{Erik Stokes}
\date{\today}

\begin{abstract}
A matroid complex is a pure complex such that every restriction is
again pure.  It is a long-standing open problem to classify all
possible $h$-vectors of such complexes.  In the case when the complex
has dimension 1 we completely resolve this question.
We also prove the
1-dimensional case of a conjecture of Stanley that all matroid
$h$-vectors are pure ${O}$-sequences.  Finally, we completely
characterize the Stanley-Reisner ideals of 1-dimensional matroid
complexes.
\end{abstract}

\maketitle
%%%%%%%%%%%%%%%%%%%%%%%%%%%%%%%%%%%%%%%%%%%%%%%%%%%%%%%%%%%%%%%%%%%%%%%
\section{Introduction}\label{sec:introduction}

Matroids are much studied objects in combinatorics, arising naturally
whenever examining independence relations between objects.  Starting
with a set of vectors in some vector space, it is natural to ask which
vectors are independent of which other vectors. This is, in some
sense, a measure of how badly our set of vectors fails to be a basis.
The collection of independent subsets of some set of vectors is knows
as a matroid.  But matroids don't only appear in connection to vector
spaces.  They can also be formed from the spanning trees of graphs,
hyperplane arrangements or any other situation where there is some
notion of ``independence'' (with the natural assumption that a subset
of an independent set is independent).  All of these objects are
matroids.  Since a subset of an independent set is independent, we can
naturally regard a matroid as a simplicial complex.

So, one way to
study matroids is by looking at the simplicial complex formed by their
independent sets (defined in Section~\ref{sec:basics}) and examining
its combinatorial, topological and algebraic (by way of the
Stanley-Reisner ideal) properties.  The complexes
that arise in this way are known as matroid complexes.  In particular,
if $\Delta$ is a simplicial complex over a finite set,
we will consider one of its fundamental invariants, the
\defword{$f$-vector}, $f(\Delta)=(f_{-1},f_0,\dotsc,f_d)$ 
where each $f_i$ is the 
number of faces of $\Delta$ with dimension $i$.  More algebraically,
we also examine the $h$-vector,
$h(\Delta)=(h_0,h_1,\dotsc,h_{d+1})$, which can be found from the Hilbert
function of the Stanley-Reisner ideal of $\Delta$.  The two sequences
are equivalent via the equations
\begin{equation}\label{eq:1}
\begin{split}
f_i&=\sum_{j=0}^{i}\binom{d-j}{i-j}h_j\\
h_i&=\sum_{j=0}^i (-1)^{i-j}\binom{d-j}{i-j} f_{i-1}.
\end{split}
\end{equation}
The algebraic study of matroid complexes was begun  by Stanley in
\cite{MR0572989} where he referred to them as ``G-complexes''.
See \cite{stanely96:_combin_and_commut_algeb} or
\cite{miller05:_combin_commut_algeb} for details on this and other
basic facts about Stanley-Reisner ideals.  The reader who is
interested in the general theory of matroids may refer to one of the
many books on the subject such as
\cite{white08:_theor_of_matroid} or \cite{oxley92:_matroid_theor}.

There are several classes of simplicial complexes whose $h$-vectors
have been completely characterized.  For example, Stanley
(\cite{MR0572989}) showed that the $h$-vectors of
shellable complexes are exactly the $O$-sequences (that is $h(\Delta)$
is the $h$-vector of an artinian monomial algebra) and likewise for the
Cohen-Macaulay complexes.  The $O$-sequences are characterized
numerically by Macaulay.  The independence complex of a matroid is
both Cohen-Macaulay and shellable so we get necessary, but not
sufficient, conditions for a given sequence to be the $h$-vector of a
matroid complex. The ultimate goal is to characterize all possible $h$-vectors
of matroid complexes and to prove the following, 30 year-old,
conjecture of Stanley \cite{stanely96:_combin_and_commut_algeb}.  

\begin{conj}\label{conj:1}
The $h$-vector of a matroid complex is a pure $O$-sequence (that is,
the $h$-vector of a artinian monomial \emph{level} algebra).
\end{conj}
This was shown for the case of ``cographic'' matroids by Merino in
\cite{MR1888777} and we prove it here for the case when $\dim\Delta=1$.

There is no known characterization of pure $O$-sequences.  However,
some necessary conditions are known beyond the requirement that the
$h$-vector be an $O$-sequence.  A theorem of Brown and Colbourn
\cite{MR1186825} states that if $h$ is the $h$-vector of a matroid
complex then 
\[(-1)^j\sum_{i=0}^j (-\alpha)^i h_i\geq 0,\;\; 0\leq j\leq d+1\]
for any $\alpha\geq 1$ with equality possible only if $\alpha=1$.
Relating to Stanley's conjecture above, Hibi (\cite{MR989204}) showed
that a pure $0$-sequence must satisfy
\[\begin{split}
h_0&\leq h_1\leq \dotsm\leq h_{\lfloor (d+1)/2\rfloor}\\
h_i&\geq h_{d-i+1},\;\; 0\leq i\leq \left\lfloor\frac{d+1}{2}\right\rfloor
\end{split}\]
and in \cite{MR1159510} conjectured that the same held for the
$h$-vectors of matroid complexes.  Also in \cite{MR1159510} Hibi proved the
weaker claim that, for the $h$-vector of a matroid,
\[h_0+h_1+\dotsm+h_i\leq h_{d+1}+h_{d}+\dotsm+h_{d+1-i}\]
for all $i\leq\lfloor \frac{d+1}{2}\rfloor$ and in
\cite{chari97:_two_decom_in_topol_combin} Chari proved Hibi's conjecture
in full.  Next, we have a result of Swartz (\cite{MR1983365})
stating that if $g_i=h_{i}-h_{i-1}$ then $g_{i+1}\leq g_{i}^{\langle
  i\rangle}$ for all $i<\tfrac{d+1}{2}$.  Another proof of this result
was given by Hausel in \cite{MR2110782}, where he also shows these
inequalities for pure $O$-sequences.
These are the best results that the author is aware of.  None of them
are sufficient for a vector to be the $h$-vector of a matroid.
For more information on what is known about pure $O$-sequences, see
\cite{boij:_shape_of_pure_o_sequen}. 

As our main results we completely classify all matroid $h$-vectors for
the case where $\dim\Delta=1$ (or equivalently, rank 2 matroids) in
Theorem~\ref{thrm:1d-binomial-h-vector}.  In fact we do somewhat more
and describe combinatorially all 1-dimensional matroid complexes
(Theorem~\ref{thrm:classify-matroids}), associating each one to a
partition of the number of vertices.  This allows us to not only
compute the $h$-vector, but also to count the number of 1-dimensional
matroid complexes on a fixed vertex set.  Using these techniques we
also obtain an algebraic description of the Stanley-Reisner ideal of
1-dimensional matroid complexes (Theorem~\ref{thrm:1d-sr-ideal}) which
we use to provide a constructive proof of Conjecture~\ref{conj:1} in
this case.

\section{Preliminary Results}\label{sec:basics}
We define a \defword{simplicial complex}, $\Delta$, over a set $X$ to
be a subset of the power set of $X$ ($2^X$) with that property that,
whenever $F\in \Delta$ and $G\subseteq F$, $G\in \Delta$.  The
elements of $\Delta$ are called \defword{faces} and the
\defword{dimension} of a face is $\dim F=\abs{F}-1$.  Faces with
dimension 0 are called \defword{vertices} and those with dimension 1
are \defword{edges}.  A $d$-face of $\Delta$ is a face with dimension
$d$ and the dimension of $\Delta$, $\dim\Delta$, is the maximum
dimension of its faces. 

 We will make use
of the following constructions on simplicial complexes.  Typically, we
take our complexes to be over the vertex set $[n]=\{1,\dotsc,n\}$.
\begin{defn}\label{def:basic-contructions}
Let $\Delta$ be a simplicial complex with vertex set $X$.  
\begin{enumerate}[(a)]
\item The \defword{$k$-skeleton} of $\Delta$ is
  $\skeleton{k}{\Delta}=\{F\in\Delta\mid\dim F\leq k\}$.
\item If $W\subs X$ then the \defword{restriction} of $\Delta$ to $W$
  is $\restrict{\Delta}{W}=\{F\in\Delta\mid F\subs
  W\}$.  If $W=X-\{v\}$ then we will write
  $\Delta_{-v}=\restrict{\Delta}{W}$ and call $\Delta_{-v}$ the
  \defword{deletion} of $\Delta$ with respect to $v$ or the deletion
  of $v$ from $\Delta$.
\item If $F\subs X$ then $\link_\Delta(F)=\{G\in\Delta\mid F\cap
  G=\emptyset, F\cup G\in\Delta\}$.  We call this the \defword{link}
  of $\Delta$ with respect to $F$.
\item If $v\not\in X$ then the cone over $\Delta$ is
  $C\Delta=\Delta\cup \{F\cup\{v\}\mid F\in\Delta\}$
\end{enumerate}
\end{defn}

We shall also make frequent use of the \defword{Stanley-Reisner} ideal
of a complex $\Delta$.  If $F\subs [n]$ then we define
$x_F=\prod_{i\in F}x_i\in S=K[x_1,\dotsc,x_n]$, for some field $K$.  The
the Stanley-Reisner ideal is the ideal
\[I_\Delta=\ideal{x_F\mid F\not\in\Delta}.\]
For information on the basic theory of Stanley-Reisner ideals we refer
the reader to \cite{miller05:_combin_commut_algeb} and
\cite{stanely96:_combin_and_commut_algeb}.  In particular, in
Section~\ref{sec:conjecture-stalely}, we will make use of the fact
that the $h$-vector $S/I_\Delta$ is the $h$-vector of $\Delta$.

We use the following,
combinatorial, definition of matroid complex.

\begin{defn}\label{def:matoid-complex}
Let $\Delta$ be a simplicial complex on $[n]$.  Then $\Delta$ is a
\defword{matroid complex} if one of the following equivalent
conditions holds.
\begin{enumerate}[(i)]
\item For every $W\subs [n]$, $\restrict{\Delta}{W}$ is pure.
\item For every $W\subs [n]$, $\restrict{\Delta}{W}$ is Cohen-Macaulay.
\item For every $W\subs [n]$, $\restrict{\Delta}{W}$ is shellable.
\end{enumerate}
\end{defn}

It is always true that (iii)$\implies$(ii)$\implies$(i).  It only remains
to show that if every subcomplex is pure then they are also
shellable.  This follows by induction, since a subcomplex of a matroid
is again matroid.  We will not do this proof here; the reader
may refer to
\cite[][Proposition~3.1]{stanely96:_combin_and_commut_algeb}.

In general, almost any construction applied to a matroid complex will
result in another matroid complex.  We summarize some of the more
useful constructions in the next elementary proposition.

\begin{prop}\label{prop:matroid-constructions}
Let $\Delta$ be a matroid complex with vertex set $[n]$.  Then the
following complexes are 
also matroid.
\begin{enumerate}[(a)]
\item $\restrict{\Delta}{W}$ for every  $W\subs [n]$
\item $C\Delta$, the cone over $\Delta$
\item $\skeleton{k}{\Delta}$, the $k$-skeleton of $\Delta$
\item $\link_\Delta(F)$ for every $F\in\Delta$.
\end{enumerate}
\end{prop}
\begin{proof}
\par\noindent
\begin{enumerate}[(a)]
\item Since
  $\restrict{(\restrict{\Delta}{W})}{V}=\restrict{\Delta}{W\cap V}$ and
  the left-hand side is, by definition, pure this follows immediately
  from the definition.
\item Let $v$ be the apex of the cone.  The cone over a pure
  complex is pure since its facets are $F\cup\{v\}$ with $F$ a facet
  of $\Delta$. So let $W\subset [n+1]$.  If $v\not\in W$ then
  $\restrict{(C\Delta)}{W}=\restrict{\Delta}{W}$, which is pure
  because $\Delta$ is matroid.  If $v\in W$ then
  $\restrict{(C\Delta)}{W}=C(\restrict{\Delta}{W}$).  By part~(a)
  $\restrict{\Delta}{W}$ is matroid and so, by induction on the number
  of vertices, $C(\restrict{\Delta}{W})$ is matroid and in particular
  pure.
\item Note that
  $\skeleton{k}{\restrict{\Delta}{W}}=\restrict{\skeleton{k}{\Delta}}{W}$.
  As in part~(b), if $W$ is a proper subset of $[n]$ then this is
  matroid, and thus pure, by induction on the number of vertices.  It
  only remains to check that $\skeleton{k}{\Delta}$ is itself pure.
  Suppose  that $\skeleton{k}{\Delta}$ has a face $F$
  with $\dim F<k$.  Since $F\in\Delta$ it must be contained in some
  facet with dimension $\dim\Delta\geq k$.  It then follows that $F$
  must be contained in some $k$-dimensional face of $\Delta$, which is
  then a face of $\skeleton{k}{\Delta}$.  Thus $F$ is not a facet of
  $\skeleton{k}{\Delta}$ and the $k$-skeleton is therefore pure.
\item This time, we check that
  $\restrict{\link_\Delta(F)}{W}=\link_{\restrict{\Delta}{W}}(F)$,
  which will then be pure by induction. We then only need to know that
  $\link_\Delta(F)$ is pure.  Suppose that $G\in\link_\Delta(F)$ is a
  facet.  Then $G\cup F\in\Delta$ must be a facet of $\Delta$.  So
  $\dim(G\cup F)=\dim\Delta$ and then $\dim G=\dim\Delta-\dim
  F-1=\dim\link_\Delta(F)$.  So the link is pure and thus matroid. 
\end{enumerate}
\end{proof}

Since if $\Delta$ has $n$ vertices and $\dim\Delta=d$ then
$\Delta_{-v}$ has $n-1$ vertices and $\link_\Delta(v)$ has dimension
at most $d-1$.  Since both the link and the deletion are smaller than
$\Delta$ is some sense, we can use induction to prove results about
$\Delta$.  Most of the arguments to follow work in the manner.

\begin{rmrk}\label{rmrk:basic-argument}
All of the statements in Proposition~\ref{prop:matroid-constructions}
follow by the same sort of argument.  First show that the desired
construction commutes with restrictions.  The proper restrictions will
then be pure by induction on the number of vertices since restrictions
of matroids are matroid.  One then only has to check that the
construction gives a pure complex.  It is important to note that the
purity of $\Delta$ does not follow from the purity of its
restrictions.  For an example, consider the complex with facets
$\{1,2\}$, $\{1,3\}$, $\{2,3,4\}$.  In this complex, all of
its proper restrictions are pure, while the complex itself is not.
\end{rmrk}

In addition to being Cohen-Macaulay, the Stanley-Reisner ideals of
matroid complexes possess  another, desirable property: they are level.
This can be see by using Hochster's Formula
\cite[][Corollary~5.12]{miller05:_combin_commut_algeb} to compute the
various degree components of the last term in the minimal free
resolution of $I_\Delta$.  Taking the link and deletion we get, by
standard results of simplicial homology, a long exact sequence.  Since
the links and deletions of matroid complexes are matroid, induction
tells us that their Stanley-Reisner ideals are level.  The long exact
sequence and Hochster's formula then forces $I_\Delta$ to be level as
well.  For a complete proof of this fact, the reader may refer to
Stanley's book \cite{stanely96:_combin_and_commut_algeb}.

In the spirit of Proposition~\ref{prop:matroid-constructions} we give
a new construction for producing new matroids from old ones.

\begin{defn}\label{def:partial-star}
Let $\Delta$ be a simplicial complex with vertex set $V$ and $v$ a
vertex of $\Delta$. 
Then we define 
\[S_v\Delta=\{E\cup\{w\}\mid E\text{ a face of }\link_\Delta(v)\}\]
where $w$ is a fixed vertex not in $\Delta$.  We consider $S_v\Delta$
as a simplicial complex with vertex set $V\cup\{w\}$.  If
$W=\{w_1,\dotsc,w_k\}$ is a set of $k$ vertices not in $\Delta$ then
\[S_v^W=S_v^k\Delta=\{E\cup\{w\}\mid E\text{ a face of }\link_\Delta(v), w\in
W\}\]
is a simplicial complex on $V\cup W$.
We call $S_v^W\Delta$ the \defword{$k$-fold partial star avoiding $v$}.
\end{defn}

\begin{figure}[ht]
\begin{minipage}[h]{.4\linewidth}
%\centering
\includegraphics[width=\linewidth]{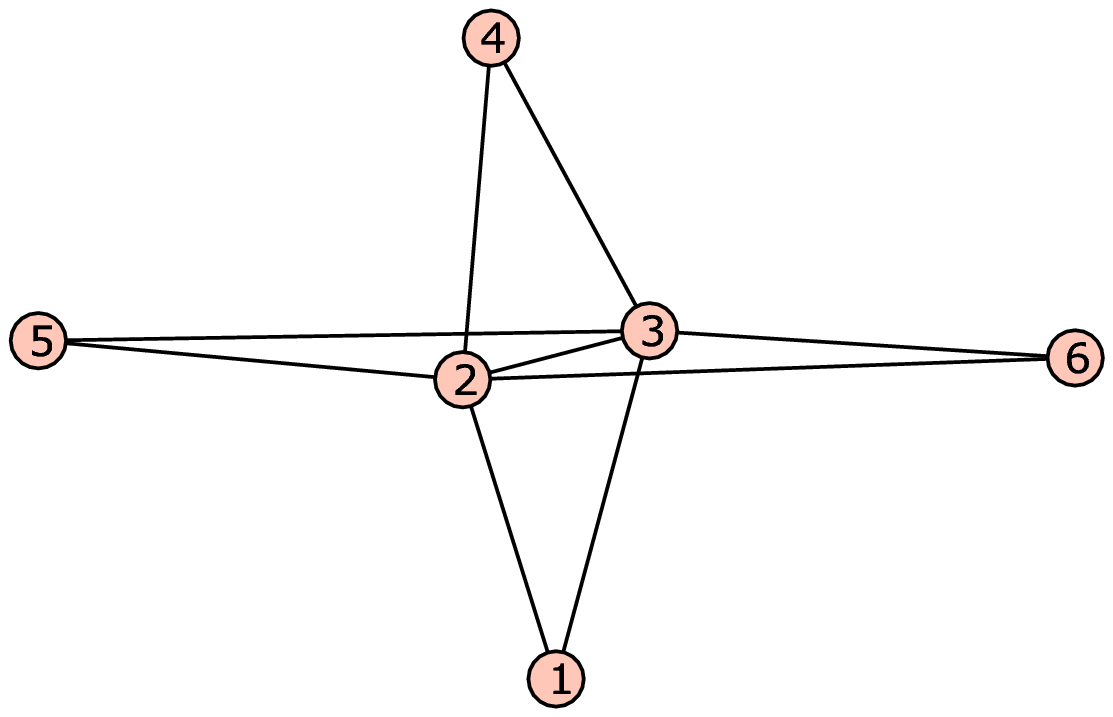}
\caption{$S^3_1K_3$}\label{fig:matroid300}
\end{minipage}
\begin{minipage}[h]{.4\linewidth}
%\centering
\includegraphics[width=\linewidth]{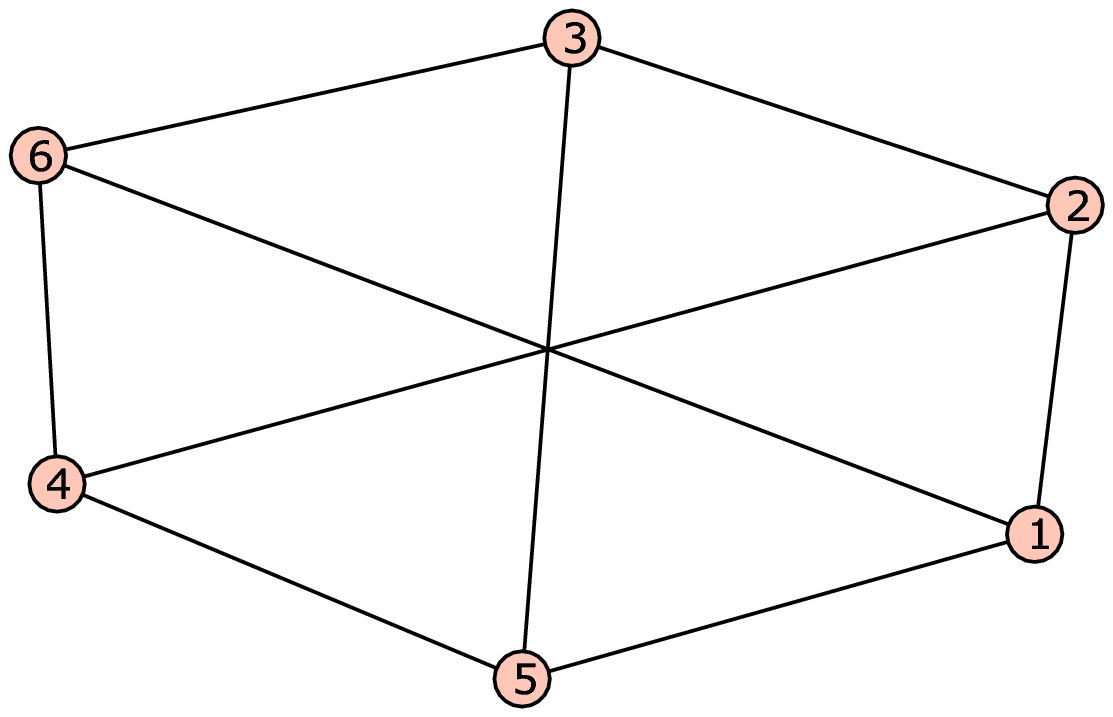}
\caption{$S^2_2 S^2_1 K_2$}\label{fig:matroid22}
\end{minipage}
\end{figure}

See Figure~\ref{fig:matroid300} for an example of the result of this
procedure.  To get this, we start with a single 3-cycle 
%(or a complete
%graph on 3 vertices if you prefer) 
and add 3 new vertices (labeled 4,5
and 6) connecting
them by an edge to every vertex except vertex 1.

Figure~\ref{fig:matroid22} represents the results of applying this
construction twice starting with a single edge between vertices 1 and
2 .  First we add vertices 3 and 4 avoiding vertex 1 and then vertices
5 and 6 avoiding vertex 2.  Note that the edges $\{3,5\}$ and
$\{3,6\}$ are in the final complex.

The next lemma informs us that, if we start with a matroid complex,
this construction will usually result in another matroid.  In fact, we
will later, in Theorem~\ref{thrm:classify-matroids}, see that all
matroids can be obtained from smaller matroids in this way.  However,
not every complex obtained by taking partial stars is matroid, even if
we start with a ``nice'' complex (for example, paths can be obtained
in this way).  So, we must impose some additional conditions on this
process.  In particular, we need to choose the vertex we avoid
properly, where the meaning of ``properly'' is given by the next
definition.

\begin{defn}\label{def:center}
Let $\Delta$ be a simplicial complex and $v$ a vertex of $\Delta$.
Then we say that $v$ is a \defword{center} of $\Delta$ if
$\link_\Delta(v)$ contains every other vertex of $\Delta$.
\end{defn}

Note that this definition depends only on the 1-skeleton of $\Delta$;
we are simply looking for vertices that are connected by an edge to
every other 
vertex.  More algebraically, if $I_\Delta$ is the Stanley-Reisner
ideal of $\Delta$, $v$ is a center of $\Delta$ if and only if  $x_v$
does not appear in any degree 2 minimal generator of $I_\Delta$.  

In Figure~\ref{fig:matroid300} vertex 2 and 3 are the only centers, while
Figure~\ref{fig:matroid22} shows us a complex that has no
centers.  If $\Delta$ is a cone then the apex of the cone is a
center.  The 
converse is false; to have a center it is only necessary that
$\skeleton{1}{\Delta}$ be the 1-skeleton of a cone.   This can be easily
seen by looking at the Stanley-Reisner ideal and using the comment in
the preceeding paragraph.

\begin{lemma}\label{lemma:matroid-attach-general}
Let $\Delta$ be a matroid complex and $W=\{w_1,\dotsc,w_k\}$ a set of vertices
not in $\Delta$.  Then $S_v^W\Delta$ is matroid if and only if $v$ is a
center of $\Delta$.
\end{lemma}
\begin{proof}
Assume that $v$ shares an edge with every other vertex of $\Delta$
(i.e., $v$ is a center of $\Delta$).  Let $\Gamma=S_v^W\Delta$.  If
$X$ is a subset of vertices of $\Gamma$ which contains only vertices
of $\Delta$ then $\restrict{\Gamma}{X}=\restrict{\Delta}{X}$, which is
pure.  So suppose that $W\cap X\neq\emptyset$.  If $v$ is the only
vertex of $\Delta$ contained in $X$ then $\restrict{\Gamma}{X}$ is 0
dimensional and thus pure.  So, assume that $X$ intersects the other
vertices of $\Delta$.  Let $X'=X\cap [\link_\Delta(v)]_0$.  Since
$\Delta$ is matroid, $\restrict{\Gamma}{X'}=\restrict{\Delta}{X'}$ is
pure.  A facet of $\restrict{\Gamma}{X}$ is of the form $\{v\}\cup E$
or $\{w_i\}\cup E$ where $E$ is a facet of $\restrict{\Gamma}{X'}$.
Since these all have the same size, $\restrict{\Gamma}{X}$ is pure and
thus $\Gamma$ is matroid.

Now, suppose that $v$ is not a center and let $a\neq v$ be a vertex of
$\Delta$ not in $\link_\Delta(v)$.  
If $\Delta$ is not pure then neither is $\Gamma=S_v^W\Delta$, so we
may as well assume that $\Delta$ is pure.
If $\Delta$ has dimension 0 then
$v$ will be a facet of $\Gamma=S_v\Delta$, which will have dimension
1.  Suppose that $\dim\Delta>1$. 
Then $\Delta_{\{a,v\}\cup W}$ has dimension 1 and has $\{v\}$ as a
facet, thus is not pure so $\Gamma$ is not matroid.
\end{proof}

The next Theorem allows us to, in many cases, reduce large matroid
complexes to much smaller ones.  This is particularly useful if in
dimension 1 where it gives a constructive procedure that produces all
matroid complexes.  

\begin{thrm}\label{thrm:classify-matroids}
Let $\Delta$ be a simplicial complex with dimension $d$.  Then
$\Delta$ is matroid if and only if $\Delta=S_{v_k}^{m_k}\dotsm
S_{v_1}^{m_1}\Gamma$ where $\Gamma$ is a matroid such that
$\skeleton{1}{\Gamma}$ is a complete graph and the $v_i$ are distinct
vertices of $\Gamma$. We allow for $\skeleton{1}{\Gamma}$ to be $K_1$,
the complete graph on 1 vertex, i.e. a point.
\end{thrm}

To prove this, we must first establish the special case when
$\dim\Delta=1$.  Then, since the skeletons of matroids are themselves
matroid, we can induct on the dimension and concern ourselves only
with the facets of $\Delta$.  That $\Delta$ has facets in the correct
places will be forced by purity.  
The next lemma is used to easily
detect the matroid-ness of a 1-dimensional complex.  This result
amounts to saying that we can walk between any 2 vertices of a
1-dimensional matroid complex by taking at most 2 steps (assuming our
steps are 1 edge long).

\begin{lemma}\label{lemma:link-facet-matroid-early}
Let $\Delta$ be a 1-dimensional simplicial complex.  Then $\Delta$ is
matroid if and only if for every vertex $v$ and every edge $E$,
$\link_\Delta(v)\cap E\neq\emptyset$.
\end{lemma}
\begin{proof}
Suppose there exists a vertex $v$ and an edge $E$ disjoint from the
link of $v$.  Let $L=[n]-\link_\Delta(v)$.  Then $\Delta_L$ has
$\{v\}$ and $E$ as facets, and so is not matroid.

Conversely, suppose that there exists a subset $W\subs [n]$ such that
$\Delta_W$ is not pure.  So $\Delta_W$ must have a 0-dimensional
facet, say $\{v\}$.  Let $v\neq w\in W$.  Since $v$ is a facet of
$\Delta_W$ we must have $\{v,w\}\not\in\Delta$.  Thus
$W\cap\link_W(v)=\emptyset$ and so any edge, $E$, of $\Delta_W$ (there
must be at least one since $\Delta_W$ is not pure) must also be
disjoint from $\link_\Delta(v)$.  Since $E$ is also an edge of
$\Delta$ the proof is complete.
\end{proof}

\begin{lemma}\label{lemma:1d-classification}
Let $\Delta$ be a matroid with dimension 1.  Then
$\Delta=S_{v_k}^{m_k}\dotsm S_{v_1}^{m_1} K_s$ where $k\leq s$ and the
$v_i$ are distinct vertices of $K_s$.
\end{lemma}
\begin{proof}
Let $v$ be a vertex of $\Delta$ and $n=f_0(\Delta)$.  We define 
$\deg v=\deg_\Delta v=\abs{\link_\Delta(v)}$.  Choose, if possible,
$v$ so that $\deg v\neq n-1$.  If there are no such vertices then
$\Delta=K_n$ and we are done.  Let $W$ be the set of vertices of
$\Delta_{-v}$ not in $\link_\Delta(v)$.  If $E\in\Delta_W$ is an edge
then $E\cap\link_\Delta(v)=\emptyset$, contradicting
Lemma~\ref{lemma:link-facet-matroid-early}.  So $\dim\Delta_W=0$.  If
$w\in W$ and $\{v,w\}\in\Delta$ then $w\in\link_\Delta(v)$, a
contradiction. Let $\Delta'$ be $\Delta$ with the vertices in $W$
deleted.  From above, we can see that any edge of $\Delta$ that is not
in $\Delta'$ must be of the form $\{w,x\}$ where $w\in W$ and
$x\in\link_\Delta(v)$. Thus, $\Delta=S_v^m\Delta'$ where $m=\abs{W}$.
Since $\Delta'$ 
has fewer vertices than $\Delta$ we can conclude by induction on $n$.
\end{proof}

Since any simplicial complex of the form $S_{v_k}^{m_k}\dotsm
S_{v_1}^{m_1} K_s$ is matroid by
Lemma~\ref{lemma:matroid-attach-general}, we now have a complete
classification of 1 dimensional matroid complexes.  Using this, we can
now prove Theorem~\ref{thrm:classify-matroids}.

\begin{proof}
Assume that $\Delta$ is matroid. We induct on $d$, the dimension of
$\Delta$.  If $d=0$ then 
$\Delta=S_v^{n-1}K_1$, where $v$ is the solitary vertex of $K_1$.  So
assume $d>0$ and choose $v\in\Delta$ to be a vertex such that
$\link_\Delta(v)$ does not contain every other vertex of $\Delta$.  If
there are no such vertices then we may set $\Delta=\Gamma$ since the
1-skeleton of $\Delta$ must be complete.  Let $W$ be the set of
vertices not in $\link_\Delta(v)$.  Let $\Delta'$ be $\Delta$ with the
vertices of $W$ deleted.  We need only show that $\Delta=S_v^W\Delta'$
since, by induction on the number of vertices, $\Delta'$ has the
required form.  By Lemma~\ref{lemma:1d-classification} we have
$[\Delta]_1=S_v^W[\Delta']_1$.  In particular $[\Delta_W]_1$ has
dimension 0.  By definition, there are no edges (and thus no higher
dimensional 
faces)  of $\Delta$ in $\{v\}\cup W$.  So, the only thing remaining to
show is that, if $E$ is a facet of $\link_\Delta(v)$ (which is
matroid and thus pure) and $w\in W$ then $\{w\}\cup E\in\Delta$.  By
induction on $d$, if $F\in\link_\Delta(v)$ is not a facet, $\{w\}\cup
F\in W$.  Suppose that $E\cup\{w\}\not\in\Delta$.  Let $X=\{v,w\}\cup
E$.  For every $e\in E$, $(E-\{e\})\cup\{w\}\in\Delta_X$ and since
$E\cup\{w\}\not\in\Delta_X$ these are all facets.  By construction,
$E\cup\{v\}\in\Delta_X$ contradicting the purity of $\Delta_X$.  Thus,
$E\cup\{w\}\in\Delta$.  Finally, we note that if $E\cup\{w\}\in W$
where $w\in W$ and $E\not\in\link_\Delta(v)$ then an identical
argument (interchanging $v$ and $w$) shows that, again,
$\Delta_{\{v,w\}\cup E}$ is not pure.  Therefore, $\Delta=S_v^W\Delta'$
and we may conclude by induction on the number of vertices.

For the converse we simply note that, by
Lemma~\ref{lemma:matroid-attach-general}, every complex of the form
$S_{v_k}^{m_k}\dotsm S_{v_1}^{m_1}\Gamma$ is matroid by
 provided that we choose the
$v_i$ so that they are centers of their respective complexes.  Being a
center depends only the 1-skeleton, which is, by assumption,
complete.  So we can simply choose the $v_i$ to be distinct vertices
of $\Gamma$ and be assured that the partial star avoiding $v_i$ is
matroid. 
\end{proof}

%%%%%%%%%%%%%%%%%%%%%%%%%%%%%%%%%%%%%%%%%%%%%%%%%%%%%%%%%%%%%%%%%%%%
\section{Dimension 1}\label{sec:dimension-1}
Our goal is to classify the h-vectors of all 1-dimensional matroid
complexes. In fact, we do something stronger and classify all
1-dimensional matroid complexes up to isomorphism in terms of
partitions.

 Since we will be working exclusively with 1-dimensional complexes, it
 will be convenient to ignore the difference between the 0-dimensional
 complex $\link_\Delta(k)$ and the set of vertices of
 $\link_\Delta(k)$.

Lemma~\ref{lemma:1d-classification} provides us with a complete
classification of 1 dimensional matroid complexes.  It only remains to
compute the possible $h$-vectors that this construction allows.  Note
that Lemma~\ref{lemma:matroid-attach-general} allows us to form a new
matroid complex $S^1_v\Delta$ whenever $\Delta$ has a center.  If
$\Delta$ has no center then we can easily give it a center by using
the next lemma.
%% We present  iterative procedures for constructing 1-dimensional
%% matroid complexes and will claim that they produce all such
%% complexes.
\begin{lemma}\label{lemma:matroid-1-cone}
Let $\Delta$ be a 1-dimensional matroid complex with h-vector
$(1,m-1,h_2)$ and $C_1\Delta$ the
1-skeleton of the cone over $\Delta$.  Then $C_1\Delta$ is matroid
with h-vector $(1,m,h_2+m)$.
\end{lemma}
\begin{proof}
Proposition~\ref{prop:matroid-constructions} shows that the cone of
$\Delta$ and its 1-skeleton, $C_1\Delta$ are matroid.
The statement about
h-vectors follows by noting that the f-vector of $C_1 \Delta$ is
$(1,m+2,f_2+m+1)$ where $f_2$ is the number of edges of $\Delta$.
\end{proof}

While the title says ``$h$-vector'', we mostly work by computing
$f$-vectors.  It will thus be convenient to use Equation~\eqref{eq:1}
to write the $h$-vector of a
complex with $f$-vector $(1,f_0,f_1)$ as  
\[h=(1,f_1-2,1-f_1+f_2).\]

So, if $\Delta$ has $h$-vector $(1,m,h_2)$ then $S_v^i\Delta$ will have
$h$-vector $(1,m+i,h_2+mi)$ since we are adding in exactly $mi$
additional edges and $i$ vertices.  If we wish to stay in the class of
matroid complexes then we must require that $\Delta$ have a center.
If it doesn't then we first apply $C_1$.
These two constructions in fact give every 1-dimensional matroid
complex on $n$ vertices and gives us a way to compute the $h$-vector.
 Here the degree of a vertex $v$ is
the number of edges containing $v$.  Note that $v$ is a center of
$\Delta$ if and only if $v$ has degree $n-1$.
\begin{thrm}\label{thrm:h-vector-of-1d-matroids}
Let $h=(1,m,h_2)$.  Then $h$ is the $h$-vector of a 1-dimensional
matroid complex if and only if one of the following holds.
\begin{enumerate}
\item\label{item:1} $h_2=x(m-x)$ for some
  $\lfloor\frac{m}{2}\rfloor\leq x\leq m$. 
\item\label{item:2} $h_2=h'+x(m-x+1)$ where
  $\lfloor\frac{m}{2}\rfloor\leq x\leq m$ 
  and $(1,x-1,h')$ is the h-vector of a matroid complex.
\end{enumerate}
\end{thrm}
\begin{proof}
Suppose that $\Delta$ is a matroid complex with $h$-vector $h$.
Let $v$ be a vertex of $\Delta$ with  degree $x+1$ and
$L=\link_\Delta(v)\cup\{v\}$. If $x=m=n-2$ and $\Delta$ is not a cone
then we may delete $v$ to obtain a matroid complex $\Delta_{-v}$ with
h-vector $(1,m-1,h')$.  So the h-vector of $\Delta$ is $(1,m,h'+m)$,
which satisfies condition~\ref{item:2} with $x=m$.  If $\Delta$ is a
cone then it satisfies condition~\ref{item:1} with $x=m$.

Assume $x\neq m$ and let $\Gamma=\restrict{\Delta}{L}$ If $w\in
L-\{v\}$ and $a\not\in L$ then $\{a,w\}\in\Delta$ by
Lemma~\ref{lemma:link-facet-matroid-early}.  So we may write
$\Delta=S_v^{[n]-L}\Gamma$.  By definition $\abs{L}=x+2$ and
$h(\Gamma_{-v})=(1,x-1,h')$ for some $h'$, as long as
$\dim\Gamma_{-v}=1$ (equivalently, as long as $\Gamma$ is not a cone).
Now simply note that $\abs{[n]-L}=n-x-2=m-x$ and that to form $\Delta$
from $\Gamma_{-v}$ we must add edges $\{a,b\}$ for every $a\in
L-\{v\}$ and $b\not\in L-\{v\}$.  There are a total of $x(m-x+1)$ such
edges.  Thus $h_2(\Delta)=h'+x(m-x+1)$.  Now suppose that $\Gamma$ is
a cone so that $h(\Gamma)=(1,x,0)$. Then from the comment just before
the proof, we see that the $h$-vector of $\Delta$ is given by
$(1,m,h_2)$ where $h_2=0+x\abs{[n]-L}=x(m-x)$.  

To get the
inequalities, we simply take $v$ to be a vertex with maximal degree.
If $x<\lfloor\tfrac{m}{2}\rfloor$ then each vertex, $w\in L$, of
$\Delta=S^{[n]-L}_v\Gamma$ has every vertex not in $L$ in its link, by
construction.  There are $m-x\geq \lfloor\tfrac{m}{2}\rfloor$ such
vertices meaning that $w$ has a larger degree than $v$.

Conversely, if $h$ satisfies one of the 2 conditions above, we must
show that there is some matroid with $h$-vector $h$.  There are
naturally 2 cases.
\begin{enumerate}[{Case} 1.]
\item If $h=x(m-x)$ then the preceding paragraph tells us how to
  construct the matroid $\Delta$.  Let $\Gamma$ be the cone over $x-1$
  vertices with apex $v$ and $\Delta=S^{m-x}_v\Gamma$.  As noted
  above, $h(\Delta)=(1,m,x(m-x))$.
\item Suppose $h=h'+x(m-x+1)$ and there is some matroid, $\Gamma$ with
  $h$-vector $(1,x-1,h')$.  Again, the needed construction is implicit
  in the preceding argument. We have that $C_1\Gamma$ is matroid with
  $h$-vector $(1,x,h'+x)$ and $\Delta=S_v^{m-x}C_1\Gamma$ has
  $x+1+1+(m-x)=m+2$ vertices and
  $h_2(\Delta)=h'+x+x(m-x)=h'+x(m-x+1)$. We choose the vertex $v$ to
  be the new vertex added when forming $C_1\Gamma$ so that we may be
  assured that it is a center and that $\Delta$ is matroid (by
  Lemma~\ref{lemma:matroid-attach-general}).
\end{enumerate}

\end{proof}

Shortly, we will give a more closed form of the same classification
(Theorem~\ref{thrm:1d-binomial-h-vector}) that
proves to be easier to work with in many cases.  However, in some
specific cases it is computationally easier to use
Theorem~\ref{thrm:h-vector-of-1d-matroids} to check that a given
sequence is a matroid $h$-vector (or not) than to use
Theorem~\ref{thrm:1d-binomial-h-vector}. 

\begin{rmrk}\label{rmrk:easy-hvectors}
Using the above theorem, we can produce several easy examples of
matroid $h$-vectors (assuming in each case that the final entry is
positive): $(1,m,m)$, $(1,m,m-1)$, $(1,m,2(m-1))$, 
$(1,m,2(m-2))$, $(1,m,3m-5)$.  The last is produced using $x=m-1$ and
$h'=m-3$, if $m\geq 3$ since $(1,m-2,m-3)$ is a matroid $h$-vector.
\end{rmrk}

\begin{rmrk}\label{rmrk:check-hvectors}
Theorem~\ref{thrm:h-vector-of-1d-matroids} provides us with a method
for checking whether or not there is a 
matroid with the specified $h$-vector.  In specific cases this can be
somewhat tedious (although it is easily automated) however, we can
eliminate certain small values immediately.
\begin{enumerate}[(i)]
\item There are no matroid $h$-vectors of the form $(1,m,h_2)$ where
  $0<h_2<m-1$ because if there were then we would also have a matroid
  $h$-vector of the form $(1,x-1,x(m-x+1))$ for some $x$.  However
  $x(m-x+1)>x(m-x)\geq m-1$ for all $1\leq x<m$ (which excludes the
  first type of $h$-vectors as well).
\item Suppose $m\geq 6$ and $m<h_2<2(m-2)$.  Then $(1,m,h_2)$ is not
  the $h$-vector of a matroid complex.  To see this, note that the
  function $g(x)=x(m-x)$ only takes on values larger than $2(m-2)$
  when $1<x<m-1$ and $g(1)=g(m-1)=m-1$, excluding $h$-vectors of the
  first type.  Similarly, the function $f(x)=x(m-x+1)$ takes on only
  values larger than $2(m-2)$ except for $f(1)=f(m)=m$.  Thus, if
  $(1,m,h_2)$ is a matroid $h$-vector then there must be another
  matroid $h$-vector $(1,m-1,h_2-m)$.  But $0<h_2-m<(m-1)-1$ and so by
  the above, there are no such $h$-vectors.
\end{enumerate}
\end{rmrk}

We now give a more closed form of
Theorem~\ref{thrm:h-vector-of-1d-matroids}.  
If $\Delta$ is a 1-dimensional matroid then, by
Lemma~\ref{lemma:1d-classification} we know that we may write $\Delta$
in the form
\[\Delta=S_{v_1}^{W_1}\dotsm S_{v_k}^{W_k}K_s.\]
From this it is straightforward to compute the $f$-vector and $h$-vector
of $\Delta$.

\begin{thrm}\label{thrm:1d-binomial-h-vector}
Let $h=(1,n-2,h_2)$, $h_2\geq 0$.  Then $h$ is the $h$-vector of a
matroid if and only if there is a sequence of numbers $m_1,m_2,\dotsc,
m_k$ such that $m_1\geq 0$, $\sum_{i=1}^k m_i=n-k$ and
\[h_2=\binom{n-1}{2}-\sum_{i=1}^k\binom{m_i+1}{2}\]
\end{thrm}
\begin{proof}
Assume $h$ is the $h$-vector of some matroid, $\Delta$. If $\Delta$ is
not the complete graph on $n$ vertices, $K_n$ (for which the claim is
obvious) then, we may write $\Delta=S_{v_1}^{W_1}\dotsm
S_{v_k}^{W_k}K_s$.  Let $m_i=\abs{W_i}$.  By construction,
$\dim\Delta_{W_i\cup\{v_i\}}=0$ and, if $X$ is not contained in any $W_i$,
$\dim\Delta_X=1$.  Moreover, all of the $W_i\cup\{v_i\}$ are pairwise disjoint.
Our construction guarantees that, if $E$ is an edge not contained in
any $W_i\cup\{v_i\}$ then $E\in\Delta$.  It follows that
$f_1(\Delta)=\binom{n}{2}-\sum\binom{m_i+1}{2}$.  It is now easy to
compute the $h$-vector of $\Delta$ and see that it is as claimed.

Conversely, if $h$ has the form given in the Theorem, we may set
$\Delta=S_{v_1}^{m_1}\dotsm S_{v_k}^{m_k}K_s$, which, as we see above,
is matroid and has the correct $h$-vector.

\end{proof}

If $m\in\N^s_0$ and $\Delta=S^{m_s}_{v_s}\dotsm S^{m_1}_{v_1} K_s$, where
the $v_i$ are distinct vertices of $K_s$ then it is easily seen that
if we choose the $v_i$ in a different order we get  isomorphic
complexes (see Lemma~\ref{lemma:permute-m-seq}).  So without loss of generality, we will always assume that
$v_i=i$ and suppress the notation.

\begin{defn}\label{def:matroid-m-seq}
If $m\in\N_0^s$ then we define $\Delta_\vec{m}=S^{m_s}\dotsm S^{m_1} K_s$
where we agree that if $m_i=0$ then $S^{m_i}\Gamma=\Gamma$ for any
complex $\Gamma$.
\end{defn}

\begin{rmrk}\label{rmrk:notation-complaint}
So, what we have (from Lemma~\ref{lemma:1d-classification}) is that
every 1-dimensional matroid is isomorphic to one of the form
$\Delta_{\vec{m}}$ for some 
sequence, $\vec{m}$, of non-negative integers with length $s$.  However, we
have chosen to construct $\Delta_{\vec{m}}$ in such a way that the
first $s$ vertices form a complete graph.  Of course, we may always
permute the vertices so that this occurs.  The point is that, while
Lemma~\ref{lemma:1d-classification} completely classifies all matroids
with dimension 1, the notation $\Delta_\vec{m}$ does not since it
implicitly assumes a particular ordering of the vertices.  Since the
author does not distinguish  between isomorphic complexes,
this  may be considered only a minor notational annoyance.

If $\vec{m}=\vec{0}$ is the zero sequence then $\Delta_\vec{0}=K_s$
and if $\vec{m}=(m_1)$ then (with the understanding that $K_1$ is a
single vertex) $\Delta_{(m_1)}$ is a 0-dimensional complex with
$m_1+1$ vertices.  Of course, all 0-dimensional complexes are matroid.
In all other cases, $\dim\Delta_{\vec{m}}=1$ as noted below in
Proposition~\ref{prop:1d-misc}~(a).
\end{rmrk}

The following is a restatement of
Theorem~\ref{thrm:1d-binomial-h-vector} using this new notation.

\begin{cor}\label{cor:1d-hvector}
If $m\in\N_0^s$ then $h$-vector of $\Delta_{\vec{m}}$ is $(1,h_1,h_2)$ where
\begin{align*}
h_1&=s+\sum_{i=1}^s m_i\\
h_2&=\binom{n-1}{2}-\sum_{i=1}^s\binom{m_i+1}{2}
\end{align*}
\end{cor}

Using an argument similar to that of
Theorem~\ref{thrm:1d-binomial-h-vector} we can describe algebraically
the structure of all possible
Stanley-Reisner ideals of 1 dimensional matroid complexes.

\begin{notation}
If $\sigma\subs [n]$ then $\m_\sigma=\ideal{x_i\mid i\in\sigma}$ (we
will write $\m=\m_{[n]}$) and $\hat{\m}_\sigma^d$ is the ideal
generated by all the squarefree monomials in $\m^d$.  If
$\abs{\sigma}<d$ then $\hat{\m}_\sigma^d=\ideal{0}$. 
\end{notation}

\begin{thrm}\label{thrm:1d-sr-ideal}
Let $I\subs S$.  Then $I$ is the Stanley-Reisner ideal of a
1-dimensional matroid on $[n]$ with $n$ vertices if and only if $I$
has the form 
\begin{equation}\label{eq:3}
I=\sum_{i=1}^k \hat{\m}_{\sigma_i}^2+\hat{\m}^3
\end{equation}
for some collection $\{\sigma_i\}$ of subsets of $[n]$ such that
$\sigma_i\cap\sigma_j=\emptyset$ whenever $i\neq j$ and
$n\geq \sum_{i=1}^k\abs{\sigma_i}$. 
\end{thrm}
\begin{proof}
Assume $I=I_\Delta$ for some 1-dimensional matroid $\Delta$ on $[n]$.
Then we may write $\Delta=\Delta_{\vec{m}}=S^{W_k}_{v_k}\dotsm
S^{W_1}_{v_1} K_s$, where $\abs{W_i}=m_i$.  We consider
$\vec{m}\in\N^s$, padding the end with 0 if needed.  Then the number
of vertices of $\Delta$ is $s+\sum m_i=n$.  Let $\sigma_i\deq
W_i\cup\{v_i\}$.  The $\sigma_i$ are all pairwise disjoint and
$\abs{\sigma_i}=m_i+1$ We then get that
$n=\sum_{i=1}^k\abs{\sigma_i}+(s-k)$, where $s-k\geq 0$.  By
construction, any edge in $\sigma_i$ is a non-face of $\Delta$ and
thus any squarefree monomial in $\hat{\m}^d_{\sigma_i}$ is in $I$.
Again the construction assures us that these are the only non-edges of
$\Delta$.  Since $\dim\Delta=1$ every degree 3 squarefree monomial
must also be in $I$.  Thus $I$ has the form given in
equation~\eqref{eq:3}.

Conversely, assume $I$ is of the form given in equation~\eqref{eq:3}.
Then set $m_i=\abs{\sigma_i}-1$ and $s=n-\sum m_i$.  Let
$m=(m_1,\dots,m_k,0,\dotsc,0)\in\N^s$.  The argument above shows us
that $I$ is the Stanley-Reisner ideal of $\Delta_{\vec{m}}$.
\end{proof}

Using our knowledge of the Stanley-Reisner ideal, we show that
$\Delta_{\vec{m}}$ is invariant up to isomorphism when the entries of
$\vec{m}$ are permuted.

\begin{lemma}\label{lemma:permute-m-seq}
Let $\tau$ be a permutation on $[s]$ and $\tau \vec{m}=(m_{\tau
  1},\dotsc,m_{\tau s})$.  Then $\Delta_{\vec{m}}\isom \Delta_{\tau\vec{ m}}$
\end{lemma}
\begin{proof}
Let $I$ and $J$ be the Stanley-Reisner ideals of $\Delta_{\vec{m}}$ and
$\Delta_{\tau \vec{m}}$ respectively.  Write
\[I=\sum_{i=1}^k \hat{\m}_{\sigma_i}^2+\hat{\m}^3\]
and
\[J=\sum_{i=1}^k \hat{\m}_{\eta_i}^2+\hat{\m}^3.\]
Clearly, $\abs{\sigma_{\tau i}}=\abs{\eta_i}$ and so since they are
all pairwise disjoint, $\abs{\cup\sigma_i}=\abs{\cup\eta_i}$.  We may
as well assume that they are equal to each other and equal to $[r]$
for some $r\leq n$.  Select a bijection $\phi_i\colon\sigma_{\tau
  i}\to\eta_i$ for each $i$.  Pasting these together (which is well
defined only because the $\sigma_i$  and $\eta_j$ are pairwise
disjoint) gives a permutation on $[n]$ (fixing everything not in
$[r]$).  This now induces an isomorphism $I\isom J$ which implies that
$\Delta_{\vec{m}}$ and $\Delta_{\tau m}$ are isomorphic as well.
\end{proof}
\begin{rmrk}\label{rmrk:h-vector-does-given-isom}
While the sequence $\vec{m}$, up to permutation, does uniquely determine the
isomorphism class of the complex $\Delta_{\vec{m}}$, the $h$-vector does not.
The matroids $\Delta_{(2,2)}$ and $\Delta_{(3,0,0)}$ both have
$h$-vector $(1,4,4)$ but are not isomorphic since $\Delta_{(3,0,0)}$
has a vertex with degree 5 but all the vertices of $\Delta_{(2,2)}$ have
degree 3.  These complexes are depicted in
Figures~\ref{fig:matroid300} and \ref{fig:matroid22}.  
\end{rmrk}

Similarly, the isomorphism class of a 1-dimensional matroid, $\Delta$,
is determined by the degree sequence of $\Delta$.  If $\Delta$ is a
1-dimensional complex (which we may regard as a graph) and $v$ is a
vertex of $\Delta$ the we define the \defword{degree} of $v$, $\deg
v=\deg_\Delta v$ to be the number of edges containing $v$, or
equivalently, the number of vertices in its link.  The \defword{degree
  sequence} of $\Delta$, $D(\Delta)$, is a sequence defined by
$D_i=\abs{\{v\in\Delta \mid \deg v=i\}}$.
\begin{lemma}\label{lemma:unique-up-to-deg-seq}
If $\Delta$ and $\Delta'$ are 1-dimensional matroids on $[n]$ then
$\Delta\isom\Delta'$ if and only if $D(\Delta)=D(\Delta')$.
\end{lemma}
\begin{proof}
It is trivial that isomorphic complexes have equal degree sequences,
so we only consider the other direction.  Since the degree sequence
determines $f_0(\Delta)$, we may assume that
$f_0(\Delta)=f_0(\Delta')=n$. Let $v$ be a vertex of $\Delta$ with
minimal degree.  If $\deg v=n-1$ then $\Delta=K_n$.  But this is the
only complex with degree sequence $D(K_n)$ ($D_{n-1}(K_n)=n$ and all
other are 0).  So assume that $\deg v<n-1$.  Let $v'$ be a vertex of
$\Delta'$ with $\deg v=\deg v'$.  Without loss of generality, we may
write $\Delta=\Delta_{\vec{m}}$.  

Since $\Delta_{\vec{m}}$ is invariant under permutations of the
entries of $\vec{m}$ we may, still without loss of generality, assume that
$v$ is in the last group of vertices to be added and likewise for $v'$
or they are the vertices being avoided.  Then $\Delta_{-v}$ and
$\Delta'_{-v'}$ have the same degree sequence (the degree of each
vertex in the link of $v$ ($v'$) goes down by 1 and the others stay
fixed).  So, by induction on the number of vertices there is an
isomorphism $\Delta_{-v}\to\Delta'_{-v'}$.  since $v$ and $v'$ must
both be ``attached'' to their deletions in the same manner (by the
construction of $\Delta_{\vec{m}}$)
this will lift to an
isomorphism simply by mapping $v\mapsto v'$.
\end{proof}
\begin{rmrk}\label{rmrk:we-can-permute}
Since we can permute the entries of $\vec{m}$ as we like, we may as
well  assume that
$m_1\geq m_2\geq\dotsm\geq m_s$.  Let $k=\max\{i\mid m_i\neq 0\}$.  If
$\Delta_{\vec{m}}$ has $n$ vertices then 
$n=s+\sum_{i=1}^s m_i=s+\sum_{i=1}^k m_i$ or equivalently,
$(m_1,\dotsc,m_k)$ is a partition of $n-s$ with length $k\leq s$.
By increasing each entry by 1, we can form a partition $\lambda=(m_1+1)+\dotsm
+(m_s+1)$ of $n$.  Thus, each 1 dimensional matroid complex
corresponds to a partition, $\lambda$, of $n$ and two matroids are
isomorphic if and only if they have the same partition.
\end{rmrk}

The partition, $\lambda$, is determined uniquely by the non-zero
entries of $\vec{m}$ and $n$.  It is often notationally more
convenient to specify $\vec{m}$ instead of $\lambda$ as $n$ is usually
understood.

\begin{ex}\label{ex:3}
Let $n=6$.  Then, as in the above remark, the partitions we are
concerned with are summarized in the following table (we allow
$\emptyset$ as the unique partition of 0).
\begin{center}
\begin{tabular}{c|c|c|c}
$n-s$&$s$&partitions & $m$\\
     &   &of $n-s$   & \\
\hline
0 & 6 &$\emptyset$&000000\\
\hline
1&5&1&10000\\
\hline
2&4&11&1100\\
&&2&2000\\
\hline 
&&111&111\\
3&3&21&210\\
&&3&300\\
\hline
&&31&31\\
4&2&22&22\\
&&4&40\\
\hline\hline
5&1&5&5
\end{tabular}
\end{center}
The last is the 0-dimensional matroid with 6 vertices.  This means we
have a total of 10 matroids on 6 vertices (see Table~\ref{tab:1}).
But there are only 8 
distinct $h$-vectors of such complexes (those that end with 10, 9, 8,
7, 6, 4, 3 and 0), such there must be either 2 $h$-vectors each with 2
matroids or a single $h$-vector with 3 matroids.  The matroids
$\Delta_{22}$ and $\Delta_{300}$ both have $h$-vector $(1,4,4)$ and
$\Delta_{2000}$ and $\Delta_{111}$ both have $h$-vector $(1,4,7)$.
See Figure~\ref{fig:111-2000}.

A similar computation with $n=7$ shows that there are a total of 14
1-dimensional matroids with 7 vertices but only 12 matroid
$h$-vectors.  As with $n=6$ there are two pairs of non-isomorphism
matroids with the same $h$-vector. In this case it is $\Delta_{3000}$
and $\Delta_{2200}$ sharing the $h$-vector $(1,5,9)$ along with
$\Delta_{1110}$ and $\Delta_{20000}$ having $h$-vector $(1,5,12)$.
\begin{figure}[h]
\begin{minipage}[h]{.4\linewidth}
\includegraphics[width=\linewidth]{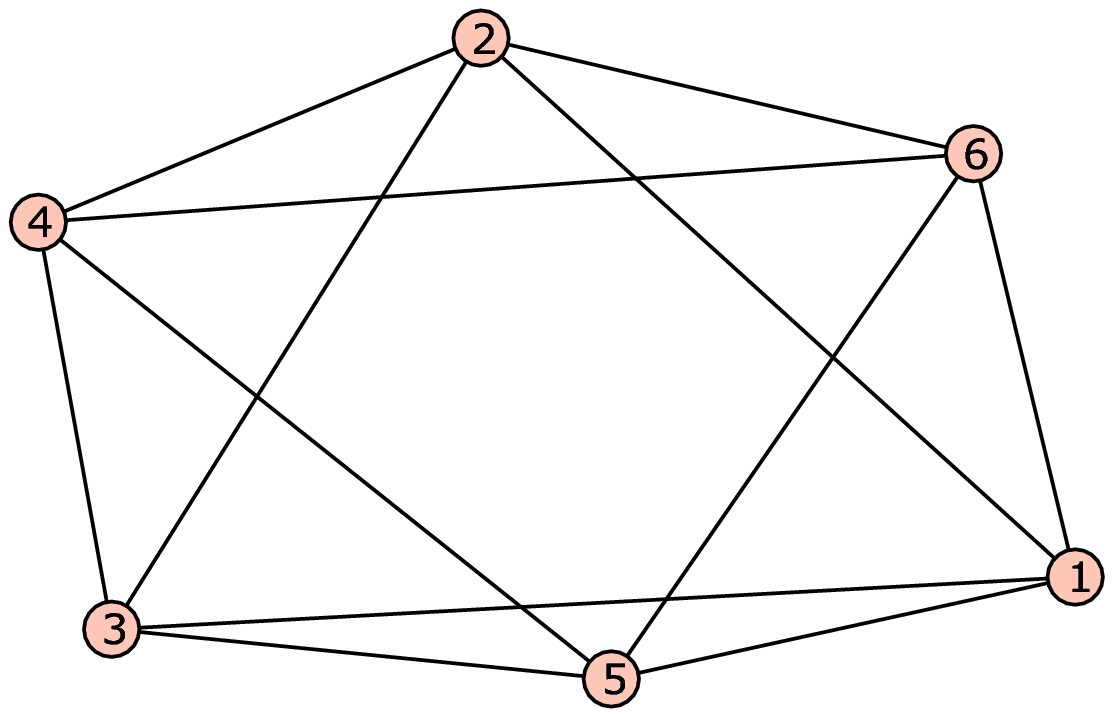}
\end{minipage}
\begin{minipage}[h]{.4\linewidth}
\includegraphics[width=\linewidth]{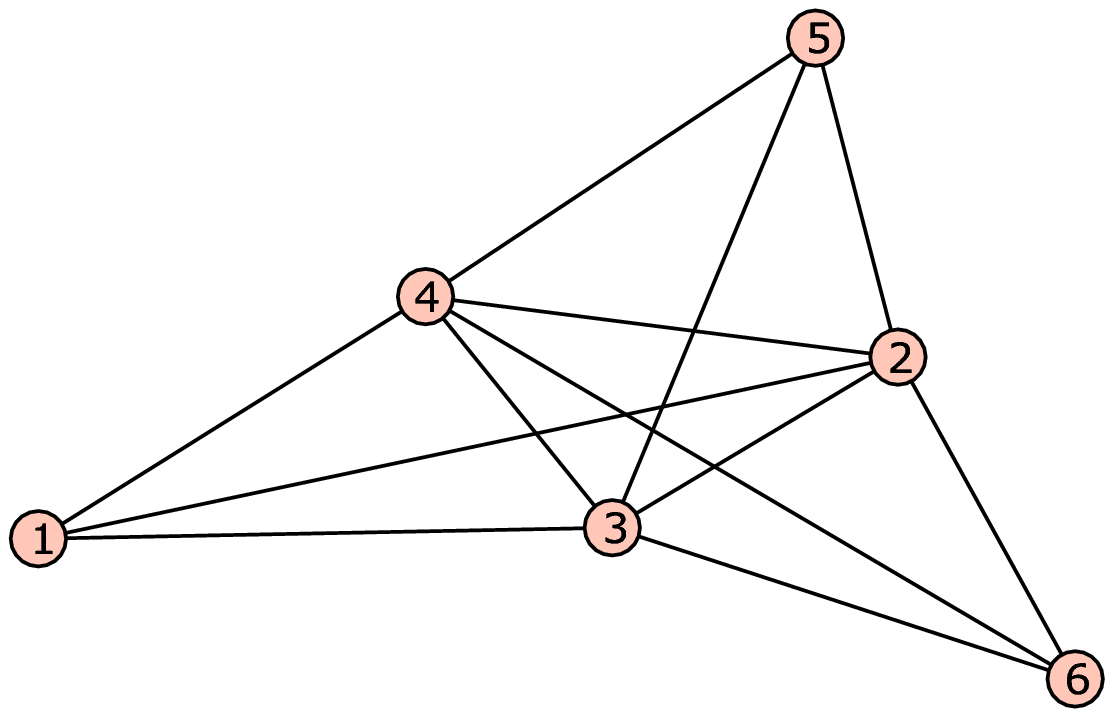}
\end{minipage}
\caption{$\Delta_{111}$ and $\Delta_{2000}$ have the same $h$-vector
  but are not isomorphic}
\label{fig:111-2000}
\end{figure}

\end{ex}

\begin{rmrk}\label{rmrk:0-dim-partition}
As we see in the above example a partition of $n$ does not necessarily
produce of 1 dimensional complex.  However the only exception is the
trivial partition $\lambda=n$, which is the complex of $n$ vertices.
We can see from the definition of $\Delta_{\vec{m}}$ that as long as the
length of $\vec{m}$ is at least 2, $\Delta_{\vec{m}}$ will contain a
complete graph on at least 2 vertices and so $\dim\Delta_{\vec{m}}=1$.
\end{rmrk}

\begin{notation}
If $\lambda$ is a partition of $n$ then we will write
$\abs{\lambda}=n$ and $\ell(\lambda)$ for the length of $\lambda$.  If
$k>1$ then $\abs[k]{\lambda}=\sum\binom{\lambda_i}{k}$  where we adopt
the convention that $\binom{a}{b}=0$ if $a<b$.
\end{notation}

\begin{defn}\label{def:partition-matroids}
\par\noindent
\begin{enumerate}[(a)]
\item If $\lambda$ is a partition of $n$ then $\Delta_\lambda$ is the
  isomorphism class of the matroid defined by the sequence
  $(\lambda_1-1,\dotsc,\lambda_{\ell(\lambda)}-1)$; $h(\lambda)$ is
  their common $h$-vector.
\item If $\Delta\isom\Delta_{\vec{m}}$ is a matroid then $\lambda_\Delta$ is the
  partition $\sum_{i=1}^s (m_i+1)$ of $n$.  We will call
  $\lambda_\Delta$ the \defword{partition associated to $\Delta$}.
\end{enumerate}
\end{defn}

We have defined $\Delta_\lambda$ so that it is not a simplicial
complex itself, but is rather a set of isomorphic  simplicial
complexes.  We will, nonetheless, continue to write things like
$\dim\Delta_\lambda$ and $h(\Delta_\lambda)$ to refer to any invariant
of the class.  We will also misuse such notation as $C\Delta_\lambda$
to refer to the class of complexes obtained from $\Delta_\lambda$ by,
in this example, taking the cone.

\begin{ex}\label{ex:1}
Consider again Figure~\ref{fig:111-2000}, which depicts the matroids
corresponding to the sequences $(1,1,1)$ and $(2,0,0,0)$ respectively.
These are elements of the classes $\Delta_{2+2+2}$ and
$\Delta_{3+1+1+1}$.  Permuting the entries of the sequences will
permute the vertices of the complexes.  This will leave $\Delta_{111}$
unchanged.  But, we can relabel the vertices so that $\{1,2,3\}$ does
not defined a complete graph, which cannot be obtained by permuting
the entries of $(1,1,1)$.  This is an element of $\Delta_{2+2+2}$ that
is not of the form $\Delta_{\vec{m}}$ for any sequence $\vec{m}$,
which does not prevent it from being matroid.
\end{ex}

\begin{thrm}\label{thrm:number-of-rank-2-matroids}
There is a bijection between isomorphism classes of  matroid complexes
with dimension at most 
1 and $n$ vertices and partitions of $n$.  In particular the number of
isomorphism classes of 1 dimensional matroids with $n$ vertices is
$p(n)-1$ where $p(n)$ is the number of partitions of $n$.
\end{thrm}
\begin{proof}
The bijection is $\lambda\mapsto\Delta_\lambda$.  The map
$\Delta\mapsto\lambda_\Delta$ is its inverse.  That every matroid with
dimension 1 can be written as $\Delta_\lambda$ is essentially the
content of Lemma~\ref{lemma:1d-classification}.
Proposition~\ref{prop:1d-misc}(a) takes care of the case when
$\dim\Delta=0$.  That these maps are inverses to each other follows
from their definitions.
\end{proof}

\begin{rmrk}\label{rmrk:if-lambda-partition}
If $\lambda\neq n$ is a partition then $h(\lambda)=(1,n-2,h_2)$ where
$n=\abs{\lambda}$ and
\[
\begin{split}
h_2(\Delta_\lambda)&=
\binom{n-1}{2}-\sum_{i=1}^{\ell(\lambda)}\binom{\lambda_i}{2}\\ 
&=\binom{n-1}{2}-\abs[2]{\lambda}.
\end{split}
\]

So two partitions, $\lambda$ and $\lambda'$ determine the same
$h$-vector if and only if $\abs{\lambda}=\abs{\lambda'}$ and
$\abs[2]{\lambda}=\abs[2]{\lambda}$.  Examples of such pairs can be
seen in Example~\ref{ex:3}.  The first such pairs are $\lambda=3+1+1$,
$\lambda'=2+2+2$ and $\lambda=3+3$, $\lambda'=4+1+1$. 

 Suppose there are $k$ entries of
$\lambda$ equal to 1 and $\lambda'$ is the partition of $n-k$ with
these entries removed.  Then
$\sum\binom{\lambda_i}{2}=\sum\binom{\lambda'_i}{2}$.  Then
$\Delta_\lambda=C_1\dotsm C_1\Delta_{\lambda'}$ (see the 
Proposition below) and $h(\lambda)$ is
easily determined from $h(\lambda')$.  So, in this sense, every
matroid $h$-vector is induced from the $h$-vector of a smaller complex,
one whose associated partition has no entry equal to 1.
\end{rmrk}

The next proposition collects various facts about the relationship
between 1-dimensional matroids and their associated partitions.  If
$\lambda$ and $\lambda'$ are partitions then we write
$\lambda+\lambda'$ for the concatenation of $\lambda$ and $\lambda'$
as a partition of $\abs{\lambda}+\abs{\lambda'}$.  Likewise,
$\vec{m},\vec{m'}$ is the concatenation of the sequences $\vec{m}$ and
$\vec{m'}$.

\begin{prop}\label{prop:1d-misc}
\par\noindent
\begin{enumerate}[(a)]
\item $\dim\Delta_\lambda=0$ if and only if  $\ell(\lambda)=1$.
\item
  $\Delta_{\vec{m},\vec{m'}}\isom[\Delta_{\vec{m}}\ast\Delta_{\vec{m'}}]_1$,
  the one-skeleton of the join. Likewise for
  $\Delta_{\lambda+\lambda'}$.
\item $\Delta_{\vec{m},0}\isom C_1\Delta_{\vec{m}}$ or equivalently
  $\Delta_{\lambda+1}=C_1\Delta_\lambda$.
\item $\Delta_\lambda$ is a cone if and only if $\lambda=(n-1)+1$
\item If $\Delta$ and $\Delta'$ are 1 dimensional matroids with $n$
  vertices then $h(\Delta)=h(\Delta')$ if and only if $\binom{n-1}{2}-
  h_2(\lambda_\Delta)= \binom{n-1}{2}-h_2(\lambda_{\Delta'})$.
   
\item\ (Klivans)\label{item:4} 
  A 1 dimensional matroid, $\Delta_{\vec{m}}$, is isomorphic to a shifted
  complex if and only if $\vec{m}$ contains at most 1 non-zero entry.
\end{enumerate}
\end{prop}
\begin{proof}
\par\noindent
\begin{enumerate}[(a)]
\item This follows immediately from the definition of
  $\Delta_\lambda$, which is formed starting with a complete graph on
  $\ell(\lambda)$ vertices.  Of course, this has dimension 0 if and
  only if $\ell(\lambda)=1$.
\item  We induct on the length of $\vec{m'}$.  If $\vec{m}$ has length 1 then
  $\dim\Delta_{\vec{m}}=0$  and, by the definition
  $\Delta_{\vec{m}}\ast\Delta_{\vec{m'}}=\Delta_{\vec{m},\vec{m'}}$.  Now, if $\vec{m'}$ has length
  $s>1$ let $\vec{m''}=(m_1,\dotsc,m_{s-1})$.  Then
  \begin{align*}
    \Delta_{\vec{m},\vec{m'}}\isom\Delta_{\vec{m},\vec{m''},m_s}&\isom[\Delta_{\vec{m},\vec{m''}}\ast\Delta_{m_s}]_1\\
    &\isom[(\Delta_{\vec{m}}\ast\Delta_{\vec{m''}})\ast\Delta_{m_s}]_1\\
    &\isom[\Delta_{\vec{m}}\ast\Delta_{\vec{m'}}]_1
  \end{align*}
  since the join is associative.
\item $C_1\Delta$ is the 1-skeleton of $C\Delta$, the join of $\Delta$
  and a single vertex.  So this follows from part (b).
\item From part (a), $\Delta_{m_1}$  has dimension 0 and from part
  (c) $\Delta_{\vec{m}}=C_1\Delta_{m_1}=C\Delta_{m_1}$  Conversely, if the
  sequence $\vec{m}$ has more than 1 non-zero entry $\Delta_{\vec{m}}$ can not be a
  cone, since it can then be written as the 1-skeleton of the join of
  two smaller complexes, at least one of which is not a cone.  So
  $\ell(\lambda_\Delta)=1$ and $\Delta_{\vec{m}}=C_1\dotsm C_1\Delta_{m_1}$,
  which is not a cone if there is more than one $C_1$.
\item This follows immediately from Corollary~\ref{cor:1d-hvector}
  after noting that $\Delta_\lambda$ has $\abs{\lambda}$
  vertices.
\item Let $\lambda=\lambda_\Delta$.  By
  Theorem~\ref{thrm:1d-sr-ideal}, we can write
  \[I=\sum_{i=1}^k \hat{\m}_{\sigma_i}^2+\hat{\m}^3\]
  where $ \abs{\sigma_i}=\lambda_i$.  We need to see that this ideal is
  squarefree strongly stable if and only if $\lambda_i=1$ for all $i>1$.  One
  direction is easy; if only $lambda_1=1$ then it is clear that
  $I_\Delta$ will be squarefree strongly stable (after permuting the
  indices that so that $\sigma_1$ is the first $\abs{\sigma_1}$
  variables).  Conversely, if $\lambda_2\neq 1$ and $x_ix_j$ is the
  product of the two variables with the smallest indices in $\sigma_2$,
  $x_{i-1}x_j\not\in I_\Delta$ since (by construction)
  $\{i-1,j\}\not\subs\sigma_2$ and it cannot be contained in any other
  $\sigma_i$ since they are all pairwise disjoint.  Thus $I_\Delta$ is
  not squarefree strongly stable whenever $\ell(\lambda)>1$ (no matter
  how we permute the indices)
\end{enumerate}
\end{proof}

\begin{rmrk}\label{rmrk:part-refit-prop}
Part~(\ref{item:4}) of Proposition~\ref{prop:1d-misc} is in fact the
same statement as Proposition~1 of \cite{klivans:_shift_matroid_compl}
which states that dimension 1 (or rank 2) shifted matroids are exactly
those obtained by starting with a dimension 0 complex and applying the
$C_1$ operator repeatedly.  This means precisely that our matroid is
isomorphism to one of the form
the form $\Delta_{(m_1,0,0,\dotsc,0)}$.  Equivalently,
$\Delta_\lambda$ contains a  shifted complex if and only if $\lambda$
has only one 
entry not equal to 1.
\end{rmrk}

If we regard the dimension 1 simplicial complex, $\Delta$, as a graph
one might want to ask about the size of the maximal cliques (that is,
maximal subsets of vertices, $W$ so that $\restrict{\Delta}{W}$ is a
complete graph).  If $\Delta$ is matroid then this is easy to determine
from looking at the partition $\lambda_\Delta$.

\begin{lemma}\label{lemma:max-cliques-from-partition-length}
If $\dim\Delta=1$ and $\Delta$ is matroid then, regarding $\Delta$ as
a graph, all the maximal cliques of $\Delta$ have
$\ell(\lambda_\Delta)$ vertices.
\end{lemma}
\begin{proof}
Let $\lambda=\lambda_\Delta$.  Of course, the sizes of the maximal
cliques depends only on the isomorphism type of $\Delta$, so we may
assume that $\Delta=\Delta_{\vec{m}}$, where $m_i=\lambda_i-1$.  Let
$s=\ell(\lambda)$.  By definition, $\restrict{\Delta}{[s]}$ is
a complete graph.  

We may assume that $\lambda$ is ordered so that
$\lambda_1\geq\lambda_2\geq\dotsm$.  If $\Delta$ is itself a complete
graph, then it is it's own unique maximal clique.  This only happens if
$\lambda=1+1+\dotsm+1$.  Ignoring this single exception, we now assume
that $\lambda_1>1$.  Let $v$ be one of the $\lambda_1-1$ vertices of
$\Delta_\vec{m}$ added in the first batch (those avoiding vertex 1).
Then
$\lambda_{\Delta_{-v}}=(\lambda_1-1)+\lambda_2+\dotsm+\lambda_s$.
This partition also has length $s$ (since $\lambda_1>1$), so all 
its maximal cliques have size $s$, by induction on the number of
vertices.  If $W$ is a maximal clique of $\Delta$ then, if $v\not\in
W$, $W$ is also a maximal clique of $\Delta_{-v}$ and therefore $\abs{W}=s$.
Suppose that $v\in W$.  Then $W-\{v\}$ is certainly a clique of
$\Delta_{-v}$, and we claim that it has size $s-1$.

To see this, we first claim that $W'=(W-\{v\})\cup\{1\}$ is a maximal
cliques and so has $s$ elements.  Let $w\in W-\{v\}$ and suppose that
$\{1,w\}\not\in\Delta_{-v}$.  Since $W$ is a clique,
$\{w,v\}\in\Delta$.  Using
  Lemma~\ref{lemma:link-facet-matroid-early}, we see that
  $\link_\Delta(1)\cup\{w,v\}$ must be in $\Delta$.  Since $\{1,w\}$
  is not in $\Delta$ we must have $\{1,v\}\in\Delta$.  But, by our
  choice of $v$, $\{1,v\}\not\in\Delta$ ($v$ is attached to $K_s$
  avoiding 1).  This contradiction indicates that our assumptions were
  wrong, so it must be that $\{1,w\}\in\Delta$, so that
  $W'=(W-\{v\})\cup\{1\}$ is a clique.

Finally, we need to show that $W'$ is a maximal clique.  Suppose that
$V$ is a maximal clique with $W'\subs V$. Then $V$ is also a clique of
$\Delta$ and $W\subs (V-\{1\})\cup\{v\}$.  Let $w\in (V-\{1\})$.  If
we can show that $\{v,w\}\in\Delta$ then $(V-\{1\})\cup\{v\}$ will be
a clique of $\Delta$ and the maximality of $W$ will imply that
$W=(V-\{1\})\cup\{v\}$, which is equivalent to $W'=V$. Since $1\in V$
and $V$ is a clique, $\{1,w\}\in\Delta$.  Then, by
Lemma~\ref{lemma:link-facet-matroid-early}, either $\{v,1\}\in\Delta$
or $\{v,w\}\in\Delta$.  Since the first is false by choice of $v$, the
second must be true.  So $W'=V$ and thus $\abs{W}=\abs{W'}=s$, as
demanded.

\end{proof}

\begin{rmrk}\label{rmrk:if-we-are}
If we are given a 1 dimensional complex, $\Delta$ and we know that it
is matroid, how can we find its associated partition,
$\lambda_\Delta$?  The answer is to search for maximal subsets
$\sigma_1,\dotsc,\sigma_s\subs [n]$ so that $\dim\Delta_{\sigma_i}=0$.  Set
$\lambda_i=\abs{\sigma_i}$ to get the associated partition. This will be
partition of $n$ since the subsets $\sigma_i$ must all be disjoint.
Why?  This is exactly the content of Theorem~\ref{thrm:1d-sr-ideal}
describing the ideal of $\Delta$ and the subset $\sigma_i$ are exactly
the subsets that appear in that theorem.

A common
problem in graph theory is to search for \defword{cliques} of a graph,
that is, subsets so that the restriction contains the maximal number
of edges.  We are doing the opposite and searching for
``anti-cliques''--- subsets whose restrictions contain the minimal
number of edges, 0.
\end{rmrk}

Within a fixed isomorphism class $\Delta_\lambda$, we can ask how many
different matroids does it contain?  This can be answered by examining
the Stanley-Reisner ideals associated to them.

Let $\lambda$ be a fixed partition of $n$.  Then we say a collection of
disjoint subsets of $[n]$, $\Omega$ is a \defword{set partition}
subordinate to $\lambda$ if $\Omega=\{\omega_1,\dotsc,\omega_k\}$
where $\abs{\omega_i}=\lambda_i$ for each $1\leq
i\leq\ell(\lambda)=k$. 

\begin{prop}\label{prop:partition-count}
Let $\lambda$ be a partition of $n$ and $s=\ell(\lambda)$.  Then there
is a bijection between the complexes in $\Delta_\lambda$ and the
collection of set partitions subordinate to $\lambda$.
\end{prop}
\begin{proof}
Let $\Delta,\Delta'\in\Delta_\lambda$.  Then, by
Theorem~\ref{thrm:1d-sr-ideal} we may write
\begin{equation}\label{eq:4}
I_\Delta=\sum_{\omega\in\Omega} \hat{\m}^2_{\omega}+\hat{\m}^3
\end{equation}
where $\Omega$ is a set partition subordinate to
$\lambda$. Conversely, any set partition defines an ideal of the same
form as \eqref{eq:4}.
The two
complexes, $\Delta$ and $\Delta'$, are \emph{equal} if and only if
$I_\Delta=I_{\Delta'}$.  But $I_\Delta=I_{\Delta'}$ if and only if they
are determined by the same set partition.
\end{proof}

\begin{rmrk}\label{rmrk:faa-di-bruno}
The number of set partitions subordinate to $\lambda$ is known to be
given by the Fa{\'a} di Bruno coefficients.  Let $a_i=\abs{\{j\mid
  \lambda_i=j\}}$ be the number of times that $i$ appears in
$\lambda$.  Then the number of set partitions subordinate to $\lambda$
is
\[\frac{n!}{a_1!a_2!\dotsm a_k! 1!^{a_1}2!^{a_2}\dotsm k!^{a_k}}.\]
A proof of this may be found in Stanley's \emph{Enumerative Combinatorics}
book \cite{stanley99:_enumer_combin} and many more details on its
history and usage in  \cite{MR1903577}.

\end{rmrk}

%%%%%%%%%%%%%%%%%%%%%%%%%%%%%%%%%%%%%%%%%%%%%%%%%%%%%%%%%%%%%%%%%%%%%%%
\section{A Conjecture of Stanley in Dimension 1}\label{sec:conjecture-stalely}
In this section, we consider the conjecture of Stanley that all matroid
$h$-vectors are pure $O$-sequences (defined below).
In dimension 1, we can
positively resolve this conjecture. Throughout this section we write
$\m=\ideal{x_1,\dotsc,x_n}$ for the maximal irrelevant ideal of $S$.
\begin{defn}\label{def:pure-ideal}
Let $I\subs S$ be a monomial ideal.  Then we say that $I$ is
\defword{pure} if the monomials outside of $I$ that are maximal with
respect to divisibility all have the same degree.  
The $h$-vectors of such
ideals are called \defword{pure $O$-sequences}.
\end{defn}

Algebraically, pure ideals are level.
%(Definition~\ref{def:linear-level}).  
An ideal is level if and only if $I\colon\m$ (the \defword{socle}) is
generated in a single degree.  This is the so-called \defword{socle
  degree}, which is the same as the twist in the last module in the
minimal free resolution (this fact is used in the proof of
Theorem~\ref{thrm:1d-stanely-conjecture}).  

\begin{lemma}\label{lemma:level-implies-pure}
Let $I\subs S$ be a pure ideal. Then $I$ is pure.
\end{lemma}
\begin{proof}
Since $I$ is level, its socle, $I\colon\m$, is generated in a single degree,
say $d$.  Let $u\not\in I$ be a monomial maximal under divisibility.
Then, by maximality $u\m\subs I$ and so $u\in I\colon\m$.  Again, by the
maximality of $u$, it must be a minimal generator of $I\colon\m$ and thus
has degree $d$.
\end{proof}

\begin{conj}[Stanley]
If $\vec{h}$ is the $h$-vector of a matroid complex then there is a pure
monomial ideal with Hilbert function $\vec{h}$.
\end{conj}

If $\Delta\in\Delta_\lambda$ is a matroid with $n$ vertices and
$\dim\Delta=0$ then we can consider the ideal $J$ generated by all
 degree 3 monomials on $x_1,\dotsc x_{n-1}$ in a polynomial
ring with $n-1$ variables.  Then $h(S/J)=(1,n-3)=h(\Delta)$.  Clearly
$J$ is artinian.  It is also strongly stable and the resolution of
such ideals is given by the Eliahou-Kervaire resolution (see
\cite[Proposition~2.12]{miller05:_combin_commut_algeb}).  This tells
us, in particular, that $J$ is level. 
%(in fact, it has a linear
%resolution).

Now, suppose that $\dim\Delta=1$ but that $\Delta$ is a cone.  Then
Proposition~\ref{prop:1d-misc} tells us that $\lambda=(n-1)+1$.  Now,
$h(\Delta)=(1,n-2,0)$, which is the $h$-vector of the ideal,
$J=\ideal{x_1,\dotsc,x_{n-2}}^3$ in a polynomial ring with $n-2$
variables.  Again $J$ has a linear resolution and is thus level.
Since these two cases are easy to handle by hand, we can from here on
out ignore them.  Note that $n$ and $(n-1)+1$ are the only partitions
of $n$ in which $n$ or $n-1$ appear.  So, we may assume that each
entry of $\lambda$ is at most $n-2$.

The following easy Lemma is a straightforward observation, but is
critical in what follows.
\begin{lemma}\label{lemma:small-seq-sum}
Suppose that $\Delta$ is a matroid with $\dim\Delta=1$.  Then $\Delta$
is the 1-skeleton of a $d$-dimensional matroid if and only if
$\ell(\lambda_\Delta)\geq d+1$.
\end{lemma}
\begin{proof}
Let $\lambda=\lambda_\Delta$ and $s=\ell(\lambda)$.  First, suppose
that $\Delta$ is the 1-skeleton of a $d$-dimensional matroid, call it
$\Gamma$.   Then $\Gamma$ contains a $d$-simplex, whose 1-skeleton is
then a complete graph on $d+1$ vertices.  Then,
Lemma~\ref{lemma:max-cliques-from-partition-length} says that
$\ell(\lambda)\geq d+1$
\end{proof}

We now prove Conjecture~\ref{conj:1} for matroids with dimension at
most 1.

\begin{thrm}\label{thrm:1d-stanely-conjecture}
Let $h$ be the $h$-vector of a matroid with dimension at most 1.  Then
there is an artinian, level monomial ideal with $h$-vector $h$ and
socle degree $n-2$ unless $\Delta$ is cone in which case it has socle
degree at most $n-3$.
\end{thrm}
\begin{proof}
Let $h=h(\Delta)$ for some matroid $\Delta$.  From the above comments,
we can assume that $\dim\Delta=1$ and that $\Delta$ is not a cone.  Let
$\lambda=\lambda_\Delta$ and
$m=(\lambda_1-1,\dots,\lambda_{\ell(\lambda)}-1)$.  From
Lemma~\ref{lemma:small-seq-sum} we know that $\ell(\lambda)\geq 2$ and
so $\sum m_i\leq n-2$.
Choose a set partition $\{\sigma_1\dotsc,\sigma_k\}$  where
$\abs{\sigma_i}=m_i$ (we ignore those $m_j$ 
equal to 0).  We can choose to do this so that $\sigma_1$ consists
of the first $m_1$ numbers, $\sigma_2$ the next $m_2$ and so on.
In this way we get a more canonical choice of set partition and we
will consider everything to depend only on the partition, $\lambda$.

Let $R=K[x_1,\dotsc,x_{n-2}]$ and $\mathfrak{n}\subs R$ be the maximal
graded ideal.  Define an ideal $J_\lambda$ by
\begin{equation}\label{eq:6}
J_\lambda=\sum \mathfrak{n}_{\sigma_i}^2+\mathfrak{n}^3\subs R.
\end{equation}
Note that this definition only make sense if $\sum m_i\leq n-2$, or
equivalently if $\dim\Delta=1$ and $\Delta$ is not a cone so that
$\cup \sigma_i\subs [n-2]$.   By
construction $J_\lambda$ is artinian so we need to show that
$h(J_\lambda)=h(\lambda)$ and that $J_\lambda$ is level.

We will work from the short exact sequence
\begin{equation}\label{eq:5}
0\rTo R/(J_\lambda:x_1)(-1)\rTo^{\cdot x_1}R/J_\lambda\rTo
R/(J_\lambda+\ideal{x_1})\rTo 0
\end{equation}
By our choice of $\{\sigma_i\}$, we can be assured that $1\in\sigma_1$
and thus that $x_1$ properly divides a minimal generator of
$J_\lambda$.  That the sequence~\eqref{eq:5} is exact is then a
standard algebraic fact. \note{prove this?}  Since both ideals on the
outside of the sequence are smaller that $J_\lambda$ we may induct to
assume that they have the proper Hilbert function.  It is thus
necessary to show that $J_\lambda+\ideal{x_1}$ and $J_\lambda:x_1$ are
of the same form as $J_\lambda$.

First, consider $J_\lambda+\ideal{x_1}$.  This is simply $J_\lambda$
with every minimal generator that $x_1$ divides removed.  That is,
\[J_\lambda+\ideal{x_1}=\sum_{i=2}^{s}\mathfrak{n}_{\sigma_i}^2
+\mathfrak{n}_{\sigma_1-\{1\}}^2+\mathfrak{n}^3\]
So $\overline{J_\lambda}=(J_\lambda+\ideal{x_1})/\ideal{x_1}\subs
\overline{R}= R/\ideal{x_1}$ is $J_{\overline{\lambda}}$ where
$\overline{\lambda}$ is the 
partition $(\lambda_1-1)+\lambda_2+\lambda_3+\dotsm$.  By induction on
 $n=\abs{\lambda}$ we get that
$h(J_{\overline{\lambda}})=h(\overline{\lambda})=h(\Delta_{-v})$ for
some vertex $v\in\Delta$.  More precisely, if $\Delta=S^{W_s}\dotsm
S^{W_1}K_s$ then we may choose $v$ to be any element of $W_1$.
But
$h(\overline{J_\lambda})=h(J_\lambda+\ideal{x_1})$ so the left side of
\eqref{eq:5} is what we need.

What about $J_\lambda\colon x_1$?  By definition, $J_\lambda$ contains
every degree 3 monomial on $\{1,\dotsc,n-2\}$, which implies that
$J_\lambda\colon x_1$ contains every degree 2 monomial on
$\{2,\dotsc,n-2\}$.  Additionally, $J_\lambda:x_1$ contains a degree 1
monomial for each element of $\sigma_1$ since $x_1x_i\in J_\lambda$
for each $i\in\sigma_1$.  So
\[J_\lambda\colon x_1=\mathfrak{n}_{\sigma_1}+\mathfrak{n}^2
=\mathfrak{n}_{\sigma_1}+\mathfrak{n}_{[n-2]-\sigma_1}^2\] and we
easily see that $J_\lambda:x_1=J_{\gamma}$ where $\gamma$ is the
partition $(n-3-\lambda_1)+1$ of $n-2-\lambda_1$.  One can see that
$h(J_\gamma)=h(\gamma)$ be computing the $h$-vectors explicitly or by
noting that $C\link_\Delta(v)$ is a cone and thus its $h$-vector,
$h(\gamma)$, is determined by the 0 dimensional complex
$\link_\Delta(v)$ and $J_\gamma$ has the proper $h$-vector by
induction on the dimension.  

So the sequence~\eqref{eq:5} tells us that
$h_i(J_\lambda)=h_{i-1}(\gamma)+h_i(\overline{\lambda})$.
Corresponding to ~\eqref{eq:5}, we have another short exact sequence
\begin{equation*}
0\rTo S/(I_\Delta\colon x_v)(-1)\rTo^{x_v}S/I_\Delta\rTo
S/I_\Delta+\ideal{x_v}\rTo 0
\end{equation*}
Since $I_\Delta\colon x_v$ is the Stalely-Reisner ideal of
$C\link_\Delta(v)$ and $I_\Delta+\ideal{x_v}$ that of $\Delta_{-v}$
which tells us that 
\[\begin{split}
h_i(\lambda)=h_i(\Delta)&=h_{i-1}(C\link_\Delta(v))+h_i(\Delta_{-v})\\
&=h_{i-1}(\gamma)+h_i(\overline{\lambda})\\
&=h_i(J_\lambda)
\end{split}
\]
and so $J_\lambda$ has the same $h$-vector as $\Delta_\lambda$.

Now, we need to show that $J_\lambda$ is level.  To do this, we apply
the mapping cone construction to sequence~\eqref{eq:5} using the fact
that, by induction, $J_\lambda:x_1$ and $J_\lambda+\ideal{x_1}$ are
level. Since all the ideals are artinian, they have projective
dimension $n-2$.  By induction on the number of variables we know the
socle degree of $J_\lambda\colon x_1$ is at most $n-3$ since
$C\link_\Delta(v)$  is a cone.  We will also know that the socle degree of
$J_\lambda+\ideal{x_1}$ is $n-2$ provided that $\Delta_{-v}$ is not a cone.
To see that this is true, recall that from
Proposition~\ref{prop:1d-misc}(d), $\Delta_{-v}$ is a cone if and only if
its associated partition  has the form $(n-2)+1$, which could happen
if $\lambda_1=n-1$ since we assumed that $\lambda_1$ was the largest
entry. This would then mean that $\Delta$ is itself a cone, a case we
excluded at the beginning.

\begin{diagram}[balance,small]
0&\rTo& R/(J_\lambda:x_1)(-1)&\rTo^{\cdot x_1}&R/J_\lambda&\rTo&
R/(J_\lambda+\ideal{x_1})&\rTo& 0\\
&&\uTo&&\uTo&&\uTo\\
&&R(-1)&\rTo&R&&R\\
&&\uTo&&\uTo&&\uTo\\
&&F_1(-1)&\rTo&G_1&&R(-1)\oplus G_1\\
&&\uTo&&\uTo&&\uTo\\
&&\vdots&&\vdots&&\vdots\\
&&\uTo&&\uTo&&\uTo\\
&&F_{n-3}(-1)&\rTo&G_{n-3}&&F_{n-4}(-1)\oplus G_{n-3}\\
&&\uTo&&\uTo&&\uTo\\
&&R(-n+2)^a&\rTo&G_{n-2}&&F_{n-3}(-1)\oplus G_{n-2}\\
&&\uTo&&\uTo&&\uTo\\
&&0&&0&&R(-n+2)^a\\
&&&&&&\uTo\\
&&&&&&0\\
\end{diagram}
We know that the final term in the rightmost column must split since
otherwise the resolution would be longer than is allowed.  Moreover,
the mapping cone construction requires that
it split with summands of $G_{n-2}$, which thus have twist
$n-2$.  No other summands of $G_{n-2}$ can split.  Any summand of
$G_{n-2}$ that does not split must also have twist $n-2$ since it is a
summand of the final term in the resolution of $J_\lambda+\ideal{x_1}$
which is $R(-n+2)^b$.  Thus every summand of $G_{n-2}$ has twist $n-2$
which means that $J_\lambda$ is level.
\end{proof}

\begin{rmrk}\label{rmrk:notice-that-proof}
Notice that in the proof of Theorem~\ref{thrm:1d-stanely-conjecture}
we do not actually need to know which $h$-vectors can occur as the
$h$-vectors of matroids.  We only need to find a class of ideals
indexed by matroids that in some sense respects links and deletions.
Since the class of matroids is closed under both operations, induction
and liberal use of the sequence~\eqref{eq:5} gets us both the
$h$-vector and levelness.  
\end{rmrk}

%%%%%%%%%%%%%%%%%%%%%%%%%%%%%%%%%%%%%%%%%%%%%%%%%%%%%%%%%%%%%%%%%%%%%%%%%%

\section{The Set of Dimension 1 Matroid
  $h$-vectors}\label{sec:set-dimension-1} 
Now we consider the collection of matroid $h$-vectors and describe
some structure this set possesses.  To begin with, we give a table
indicating which, out of all Cohen-Macaulay $h$-vectors, are matroid.
The $h$-vector $(1,n-2,h_2)$ is a Cohen-Macaulay $h$-vector if and
only if $h_2\geq 0$.  In fact, a 1-dimensional simplicial complex is
Cohen-Macaulay if and only if it is connected, which is true if and
only if $h_2\geq 0$.

In Table~\ref{tab:1}, the row number of each entry corresponds to the
number of variables.  The possible Cohen-Macaulay $h$-vectors are
listed with the maximal values ($\binom{n-1}{2}$) being aligned on the
left side.  Those entries that are $h_2$ for a matroid with $n$
vertices are shaded.  The unshaded entries are not matroid $h$-vectors.

\begin{sidewaystable}
\centering
\caption{Matroid $h$-vectors}\label{tab:1}
\small{
\begin{tabular}{c|ccccccccccccccccccccccccccccc}
$n$&$h_2$\\
\hline\\
2&\textbf{0}\\
3&\textbf{1}&\textbf{0}\\
4&\textbf{3}&\textbf{2}&\textbf{1}&\textbf{0}\\
5&\textbf{6}&\textbf{5}&\textbf{4}&\textbf{3}&\textbf{2}&1&\textbf{0}\\
6&\textbf{10}&\textbf{9}&\textbf{8}&\textbf{7}&\textbf{6}&5&\textbf{4}&\textbf{3}&2&1&\textbf{0}\\
7&\textbf{15}&\textbf{14}&\textbf{13}&\textbf{12}&\textbf{11}&\textbf{10}&\textbf{9}&\textbf{8}&7&\textbf{6}&\textbf{5}&\textbf{4}&3&2&1&\textbf{0}\\
8&\textbf{21}&\textbf{20}&\textbf{19}&\textbf{18}&\textbf{17}&\textbf{16}&\textbf{15}&\textbf{14}&\textbf{13}&\textbf{12}&\textbf{11}&\textbf{10}&\textbf{9}&\textbf{8}&7&\textbf{6}&\textbf{5}&4&3&2&1&\textbf{0}\\
9&\textbf{28}&\textbf{27}&\textbf{26}&\textbf{25}&\textbf{24}&\textbf{23}&\textbf{22}&\textbf{21}&\textbf{20}&\textbf{19}&\textbf{18}&\textbf{17}&\textbf{16}&\textbf{15}&14&\textbf{13}&\textbf{12}&11&\textbf{10}&9&8&\textbf{7}&\textbf{6}&5&4&3&2&1&\textbf{0}\\
\end{tabular}

}
\end{sidewaystable}

The first 2 rows of this table are automatic: there is only a single
1-dimensional complex with 2 vertices and only 2 pure complexes with 3
vertices.  All of these are matroid and have $h$-vectors $(1,0,0)$,
$(1,1,0)$ and $(1,1,1)$.  Moreover, for any $n$, there is a matroid
with $h$-vector $(1,n-2,0)$, namely, the cone over $n-1$ vertices.  So
we may shade the 0 on each row as well.  From these, one may
completely fill in the rest of the table.  For notational convenience,
we will write $m=n-2$.

Recall that if $\dim\Delta=1$ then $C_1\Delta$ is the 1-skeleton of
the cone over $\Delta$.  By Lemma~\ref{lemma:matroid-1-cone}, whenever
$\Delta$ is matroid so is $C_1\Delta$ and
$h(\Delta)=(1,m,h_2)$ then $h(C_1\Delta)=(1,m+1,h_2+m+1)$.  Writing
$h_2=\binom{n-1}{2}-k$ for some $k$ we see that
$h_2+m=\binom{n}{2}-k$.  This is the entry in Table~\ref{tab:1}
directly below that of $h(\Delta)$.  So if we have a shaded entry in
Table~\ref{tab:1}, we may also shade each entry directly below it.

This still gives only a small portion of Table~\ref{tab:1}.  To fill
in the rest, we need another operation.  Fortunately, we have one.
Recall the definition of the partial star, $S_v^k$
(Definition~\ref{def:partial-star}). If $\Delta$ is any matroid with a
center (equivalently, $\Delta=C_1\Gamma$ for some other complex
$\Gamma$) then $S_v^k\Delta$ is again matroid.  We earlier computed
the $h$-vector of $S_v^k\Delta$.  If $k=1$ then this gives us a
``move'' from $(1,m,h_2)$ to $(1,m+1,h_2+m)$, which lies diagonally
down and to the right.  If $k=2$ we first move down one and to the
right one step and then down one step and to the right 2.  Continue,
moving an additional step to the right each time.  This is illustrated
below; the $\times$ indicates a matroid $h$-vector.  Note that, since
$\Delta$ must contain a center, we can only 
begin with an $h$-vector that has another $h$-vector directly above it,
or with a 0.

\begin{tabular}{ccccccccc}
$\times$&&&\\
$\times$&&&\\
--&$\times$&&\\
--&--&--&$\times$\\
--&--&--&--&--&--&$\times$\\
\end{tabular}

Thus, we may move straight down, or in parabolic arcs running parallel
to the upper edge of Table~\ref{tab:1}.  Everything we hit is
guaranteed to be matroid by the results of the previous sections.
That this gives all matroid $h$-vectors is the content of our
classification of 1-dimensional matroids.  Let's fill in the first few
rows as an example.

\begin{ex}\label{ex:11}
We will fill in the first 6 rows of the Table~\ref{tab:1}.  We begin
with all entries unshaded.
\begin{center}
\begin{tabular}{c|ccccccccccccccccccccccccccccc}
$n$&$h_2$\\
\hline\\
2&{0}\\
3&{1}&{0}\\
4&{3}&{2}&{1}&{0}\\
5&{6}&{5}&{4}&{3}&{2}&1&{0}\\
6&{10}&{9}&{8}&{7}&{6}&5&{4}&{3}&2&1&{0}\\
\end{tabular}
\end{center}
We obtain the first row for free since the only
1-dimensional complex with 2 vertices is matroid.  So, we first shade
in the 0 in the first row.  In fact, we may shade the 0 in each row,
since they all are matroid $h$-vectors (as they are all cones).

\begin{center}
\begin{tabular}{c|ccccccccccccccccccccccccccccc}
$n$&$h_2$\\
\hline\\
2&\textbf{0}\\
3&{1}&\textbf{0}\\
4&{3}&{2}&{1}&\textbf{0}\\
5&{6}&{5}&{4}&{3}&{2}&1&\textbf{0}\\
6&{10}&{9}&{8}&{7}&{6}&5&{4}&{3}&2&1&\textbf{0}\\
\end{tabular}
\end{center}

As in the discussion above, we may move directly down from any matroid
$h$-vector and get another matroid $h$-vector.  This gives us some
additional shaded entries, which we put a box around to distinguish
from those obtained in the previous step.

\begin{center}
\begin{tabular}{c|ccccccccccccccccccccccccccccc}
$n$&$h_2$\\
\hline\\
2&\textbf{0}\\
3&\fbox{\textbf{1}}&\textbf{0}\\
4&\fbox{\textbf{3}}&\fbox{\textbf{2}}&{1}&\textbf{0}\\
5&\fbox{\textbf{6}}&\fbox{\textbf{5}}&{4}&\fbox{\textbf{3}}&{2}&1&\textbf{0}\\
6&\fbox{\textbf{10}}&\fbox{\textbf{9}}&{8}&\fbox{\textbf{7}}&{6}&5&\fbox{\textbf{4}}&{3}&2&1&\textbf{0}\\
\end{tabular}
\end{center}

Now we begin to move diagonally.  The first entry at which we may
begin this is the 1 on the second row.  However, this gives us no new
entries.  The next choice is the 0 on the second row.  From here we
get more new entries.

\begin{center}
\begin{tabular}{c|ccccccccccccccccccccccccccccc}
$n$&$h_2$\\
\hline\\
2&\textbf{0}\\
3&{\textbf{1}}&\textbf{0}\\
4&{\textbf{3}}&{\textbf{2}}&\fbox{\textbf{1}}&\textbf{0}\\
5&{\textbf{6}}&{\textbf{5}}&{4}&{\textbf{3}}&\fbox{\textbf{2}}&1&\textbf{0}\\
6&{\textbf{10}}&{\textbf{9}}&{8}&{\textbf{7}}&{6}&5&{\textbf{4}}&\fbox{\textbf{3}}&2&1&\textbf{0}\\
\end{tabular}
\end{center}

We do this again with another entry.  We can not use the 1 on the
third row as it does not have a shaded entry above it.  We may,
however use the 2 on the third row.

\begin{center}
\begin{tabular}{c|ccccccccccccccccccccccccccccc}
$n$&$h_2$\\
\hline\\
2&\textbf{0}\\
3&{\textbf{1}}&\textbf{0}\\
4&{\textbf{3}}&{\textbf{2}}&{\textbf{1}}&\textbf{0}\\
5&{\textbf{6}}&{\textbf{5}}&\fbox{\textbf{4}}&{\textbf{3}}&{\textbf{2}}&1&\textbf{0}\\
6&{\textbf{10}}&{\textbf{9}}&{8}&{\textbf{7}}&\fbox{\textbf{6}}&5&{\textbf{4}}&{\textbf{3}}&2&1&\textbf{0}\\
\end{tabular}
\end{center}
 We may shade one more entry, the 8 on the last row as it lies
 directly below a matroid.  This completes the table as every other
 valid move will land on an already shaded entry.

\begin{center}
\begin{tabular}{c|ccccccccccccccccccccccccccccc}
$n$&$h_2$\\
\hline\\
2&\textbf{0}\\
3&{\textbf{1}}&\textbf{0}\\
4&{\textbf{3}}&{\textbf{2}}&{\textbf{1}}&\textbf{0}\\
5&{\textbf{6}}&{\textbf{5}}&{\textbf{4}}&{\textbf{3}}&{\textbf{2}}&1&\textbf{0}\\
6&{\textbf{10}}&{\textbf{9}}&\fbox{\textbf{8}}&{\textbf{7}}&{\textbf{6}}&5&{\textbf{4}}&{\textbf{3}}&2&1&\textbf{0}\\
\end{tabular}
\end{center}
\end{ex}

Notice in Example~\ref{ex:11} that many $h$-vectors can be reach in
multiple ways.  For example $(1,3,4)$ lies below $(1,2,1)$ and
diagonally from $(1,2,2)$.  It is often, but not always, true that
different paths result in non-isomorphic matroids.  In this case,
there is only one matroid with $h$-vector $(1,3,4)$.

How are the moves described above reflected in the associated
partitions?  From Proposition~\ref{prop:1d-misc} we see that a move
directly down simply adds a 1 to the end of the partition.  The
diagonal moves are move complex.  We may only apply these to
partitions that contain a 1 (we will assume it is written last).
Then, each time we move down a row, this 1 is increased by 1.  That is
we move from the partition $3+1+1$ to $3+1+2$ and then to $3+1+3$ and
so on.  At times, there will be multiple partitions in a given space.
This will occur if and only if the matroids they produce have the same
$h$-vector.

Below, we give a table indicating where one Table~\ref{tab:1} the
associated partitions are located as well as the number of matroid
complexes with a specified $h$-vector.  As a space saving measure, we
will  not use the $+$ between the terms of the partitions and we will
not list more than a single 1. We will use a subscript to indicate the
number of times a value is repeated.  For example, $321_2=3+2+1+1$.

\begin{table}[h!]
\centering
\caption{Associated partitions sorted by $h$-vector}\label{tab:2}
\begin{tabular}{c|lllllllllll}
$n$&$h_2$\\
\hline\\
2&$1_2$\\
3&$1_3$&21\\
4&$1_4$&$21_2$&22&31\\
5&$1_5$&$21_3$&221&$31_2$&32&--&41\\
6&$1_6$&$21_4$&$221_2$&$\stackrel{31_3}{\scriptstyle{222}}$&321&--&$\stackrel{41_2}{\scriptstyle{33}}$&42&--&--&51 
\end{tabular}
\end{table}

Note in Table~\ref{tab:2} that there are only two spaces (and so only
two corresponding matroid $h$-vector) containing more than one
partition. This matches what we have seen before in Example~\ref{ex:3}
that there are only two matroid $h$-vectors that are the $h$-vectors of
two different complexes.  If we were to continue Table~\ref{tab:2} to
row 7 then we would get another pair of partitions with the same
$h$-vectors, 2221 and $31_4$ together with $41_3$ and 331.  These lie
directly below the duplicated pairs on row number 6.  This will occur
again with 8 vertices and we also get a pair by moving diagonally
since on row 7 there is a space all of whose partitions contain a 1.
We also get another new pair $321_3$ and $2222$ from where a diagonal
move and a vertical move happen to coincide.

%%%%%%%%%%%%%%%%%%%%%%%%%%%%%%%%%%%%%%%%%%%%%%%%%%%%%%%%%%%%%%%%%%%%%%%%

\def\cprime{$'$}
\providecommand{\bysame}{\leavevmode\hbox to3em{\hrulefill}\thinspace}
\providecommand{\MR}{\relax\ifhmode\unskip\space\fi MR }
\providecommand{\MRhref}[2]{%
  \href{http://www.ams.org/mathscinet-getitem?mr=#1}{#2}
}
\providecommand{\href}[2]{#2}

%%%%%%%%%%%%%%%%%%%%%%%%%%%%%%%%%%%%%%%%%%%%%%%%%%%%%%%%%%%%%%%%%%%%%%%
\end{document}